\newcommand{\ext}[1]{\stackrel{#1}{\wedge}}

\def\QQ{\rlap {\raise 0.4ex \hbox{$\scriptscriptstyle |$}}
  {\hskip -0.1em Q}}

\def\CC{{\rm \kern.24em \vrule width.02em height1.4ex
depth-.05ex \kern-.26em C}}
\def\C*{{\CC}^{*}}
\def\CQ{{{\CC}^{*}} \otimes {\QQ}}
\def\Z{{\bf Z}}

\def\ra{\rightarrow }
\def\ify{\infty}

\def\A{{\cal A}}

\def\G{{\cal G}}
\def\L{{\cal L}}

\def\P{{\cal P}}
\def\T{{\cal T}}
\def\X{{\cal X}}

\def\Tg{{\cal T}_g}
\def\Th{{\cal T}_h}
\def\Tt{{\cal T}_{\tilde g}}
\def\Tin{{\cal T}_{\infty}}
\def\TinX{{\cal T}_{\infty}(X)}
\def\TinY{{\cal T}_{\infty}(Y)}
\def\Tpi{{\cal T}(\pi)}

\def\Ch{{\cal C}_h}

\def\Hin{{H}_{\infty}}
\def\THin{{\cal T}({H}_{\infty})}
\def\THinX{{\cal T}({H}_{\infty}(X))}

\def\Mg{{\cal M}_g}
\def\Mh{{\cal M}_h}

\def\Met{{Met(\Sigma)}}
\def\MCh{{MC}_{h}}
\def\MCin{{MC}_{\infty}}

\def\Dng{{DET}_{n,g}}
\def\Dnt{{DET}_{n,\tilde g}}

\def\D0Q{{DET}(0,{\QQ})}

\def\tg{{\tilde g}}
\def\tX{{\tilde X}}

\def\a{\alpha }
\def\b{\beta }
\def\ga{\gamma }
\def\Ga{\Gamma }
\def\d{\delta }
\def\r{\rho }
\def\p{\psi }

\def\la{\lambda }

\def\s{\sigma }
\def\S{\Sigma }

\def\om{\omega }
\def\Om{\Omega }

\long\def\comment#1\endcomment{}

\def\mapright#1{\smash{
    \mathop{\longrightarrow}\limits^{#1}}}

\def\mapdown#1{\Big\downarrow
 \rlap{$\vcenter{\hbox{$\scriptstyle#1$}}$}}
\documentstyle[12pt]{article}
\begin{document}
\baselineskip16pt

\begin{center}
{\bf MATHEMATICS IN AND OUT OF}\\
{\bf STRING THEORY}
\footnote{Expanded version of the opening lecture in the 37th International
Taniguchi Symposium: {\it ``Topology and Teichm\"uller spaces''}, July 1995.
{}~~~~~~~~~~~~~~~~~~~~~~~~~~~~~~~~~~~~~~~~~~[Preprint no.:{\bf imsc/95/32}]}\\

\vspace{.35cm}

{\bf Subhashis Nag}\\
{The Institute of Mathematical Sciences} \\
{Madras 600113, INDIA}\\

\vspace{.6cm}
**********************************************\\
\vspace{-.16cm}
{\em Dedicated to the memory of my father and mother}\\
***********************************************\\
\end{center}

\vspace{.3cm}

The upsurge of excitement amongst theoretical physicists,
over the subject of string theory, has filtered through to an
appreciable extent into the mathematical community. Whereas
the basic reason for the excitement amongst the physicists has been the
vision of unification for the fundamental forces of nature via string
theory, in mathematics its interest has been the wide range of deep
ideas that have been involved, and the fascinating interconnections
that have emerged.

``Stringy'' ideas in mathematics include :

\noindent
{\bf (a)} The investigation of path integrals over the spaces of Riemann
surfaces, leading to a natural modular-invariant measure (``the
Polyakov volume form'') on the Teichm\"uller spaces. The description of
this measure by complex geometry of the Teichm\"uller spaces, involving
the Mumford isomorphisms.

\noindent
{\bf (b)} A search for infinite-dimensional
``Universal Teichm\"uller spaces'' of Riemann surfaces that parametrize
simultaneously complex structures on surfaces of all topologies,
-- and the canonical relationships between various natural candidates
for such a moduli space. This is important for a non-perturbative
formulation of string theory.

\noindent
{\bf (c)} The study of the unitary (and projective unitary)
representations of the diffeomorphism group of the circle (the closed
string!); at the infinitesimal level, this is the representation
theory of the Virasoro algebra. Indeed, there is an intimate
relationship ([NV]) between the group $Diff(S^1)$ and the Teichm\"uller
spaces -- demonstrating that (c) is deeply related to (b).

It goes without saying that the above topics by no means exhaust the
mathematical challenges raised by string theory. Owing to restrictions
of space and time, and, more importantly, of the author's knowledge,
we shall deal in these notes only with some matters pertaining to items
(a) and (b) above. For more directions, see the references cited.
We will provide here an exposition of the Polyakov-Mumford
construction on the Teichm\"uller space $\Th$
of Riemann surfaces of fixed genus $h \geq 2$; we will then explain
our recent work (with Indranil Biswas and Dennis Sullivan) coherently
fitting together this construction over the Universal Commensurability
Teichm\"uller space, $\Tin$. In fact, $\Tin$ qualifies as a parameter space
of the type desired in (b), because it comprises compact Riemann surfaces
of all genus.

The connecting thread intertwining all the mathematics we discuss is
the natural appearance of the Teichm\"uller/moduli spaces of Riemann
surfaces. The most fundamental point is that, mathematically
speaking, the quantum theory of the closed bosonic string is the theory
of a sum over random surfaces ("world-sheets") that are swept out by
strings propagating in spacetime. Owing to conformal invariance
properties of the ``Polyakov action'', that sum finally reduces to an
integral over the parameter space of Riemann surfaces. {\it Namely,
the quantum string theory, as a sum over random surfaces, precipitates
a natural measure -- the Polyakov measure -- on each moduli space
$\Mh$.}

We adopt the attitude that we are addressing mathematicians with
no prior exposure to (the pulling of) strings. We have taken particular
pains (Sections I and II) to explain to a mathematical audience the
reduction of the Polyakov functional integral from an ill-defined and
infinite dimensional situation to a finite dimensional one.

\noindent {\bf Acknowledgements:} It is a pleasure to record my
thanks to the {\it Taniguchi Foundation of Japan} for inviting me to speak
at the 37th Taniguchi Mathematics Symposium: ``Topology and Teichm\"uller
Spaces''  in Finland, and for inviting this exposition for the Proceedings.
The generous financial support given by the Foundation provides me
very pleasant memories of my travels in Finland during July-August, 1995.
I thank sincerely all the many mathematicians who have listened to my
talks over the years and have helped me to understand the material
being presented here.

The paper is organized as follows:

\begin{tabular}{lcl}
Section I: Polyakov action and the string path integral. \\
Section II: Polyakov volume form on the moduli space.\\
Section III: Geometry of Teichm\"uller space and Polyakov volume.\\
Section IV: The Universal Teichm\" uller space of compact surfaces.\\
Section V: Universal Polyakov-Mumford on $\Tin$. \\
\end{tabular}

\vspace{.3cm}
\noindent
{\it {\bf I:} POLYAKOV ACTION AND THE STRING PATH INTEGRAL}

\noindent
{\bf I.A. Mechanics:}
The uninitiated mathematician may not object to being reminded of how
path integrals arise in the first place. Mechanics, whether classical
or quantum can be formulated as arising from a Lagrangian that
allows one to assign a weight, called the ``action'', to any choice of
an admissible path in configuration space connecting given initial and
final boundary conditions.

For instance, suppose a particle (or a mechanical system) is moving
in configuration space ${\bf R}^d$, with position at time $t$ being
$x(t) = (x_1(t), \ldots , x_d(t))$ from $x(t_o) = x_o$ to $x(t_1)
= x_1$. Then the Lagrangian, $L(x(t), \dot{x}(t), t)$, is a functional
of the path and the {\it action} for that choice of path
is defined by :
$$
S(x(t)) = \int^{t_1}_{t_0} L(x(t), \dot{x}(t), t) dt \eqno(1.1)
$$

\noindent
{\bf Example:} Particle moving in a potential
$V:{\bf R}^d \rightarrow {\bf R}$. $L$ may be taken as
``kinetic energy minus potential energy'', namely,
$L = \frac{m}{2} \| \dot{x}(t) \|^2 - V(x(t))$.

Amongst all admissible paths interpolating between $x(t_0) =
x_0$ and $x(t_1) = x_1$, the actual classical path followed
by the particle is the ``best'' path - namely the value of the
action should be extremal (minimal). Thus the classical equations
of motion are the Euler-Lagrange equations for the variational
problem of minimising the action (1.1). Applied to the case of
the example above, the reader can easily check that the
Euler-Lagrange equations are simply Newton's equations of motion.

In quantum mechanics [Feynman's path integral formulation] there
is no determined ``preferred'' path from $(x_0, t_0)$ to
$(x_1,t_1)$. Rather ``all'' joining paths are possible histories
of transition, and the proper question therefore is not which
path the particle follows, but what is the probability amplitude
that a particle at $x_0$ (at time $t_1$) will be at $x_1$
(at time $t_1$). That probability is taken to be a certain
weighted average over all interpolating paths, where a path
is weighted by $exp(-Action(path))$. Thus notice that the classical
path (with the minimum action) is given the highest weight, and
paths near to the classical one would get relatively high weighting.
The {\it path integral} answering the basic quantum query is
$$
Z = \int_{\{x(t)\}} e^{-S(x(t))} Dx \eqno(1.2)
$$

\noindent where $\{x(t)\}$ is the family of all admissible (say
continuous) paths joining the given initial and final conditions,
and $Dx$ represents some Wiener-type measure on this path-family.

The analog of this infinite-dimensional path integral is what we
will now describe for Polyakov's bosonic string. The crucial discovery
is that in a particularly happy situation (namely when the spacetime
dimension $d=26$), the integral reduces from the infinite-dimensional
space of possible paths (=world-sheets) to the finite-dimensional
moduli spaces of Riemann surfaces.

\vspace{.3cm}
\noindent {\bf I.B. String Theory :}

{\it String theory is the theory of fundamental particles considered
as being one dimensional, ``strings'', rather than as zero dimensional
(``point-like'') objects}. Thus a closed string, (therefore a {\it circle}
-- there being only one closed 1-manifold), propagating in a spacetime
${\bf R}^d$ sweeps out a 2-dimensional surface called its
{\it world-sheet}. Several strings can interact, and more may be
created or annihilated, without the world-sheet (which is a
2-manifold) becoming singular.

\begin{center}
\vspace{2cm}

{\bf $\Sigma$}

\vspace{2.2cm}

{{\it Figure I.1:} {\bf String world-sheet}}

\end{center}

The figure illustrates that a configuration of two strigns at an
initial time can become (say) three at a later time while sweeping
out a {\it non-singular} history $\Sigma$. Compare this with the case
for point-like particles!

In Polyakov's string theory [Pol1] the action assigned to any
particular world-sheet $\Sigma$ depends on its location (embedding)
in ${\bf R}^d$ as well as on an arbitrarily chosen smooth Riemannian
metric on $\Sigma$. Both these freedoms, in the choice of embedding
and metric, have to be integrated out (i.e., averaged over) in
setting up the path integral.

The fundamental path integral that one needs to evaluate
is the ``vacuum to vacuum'' amplitude -- meaning that both the
initial and final configurations are taken to be empty. Consequently,
the world-sheets are without boundary -- namely {\it closed surfaces}
embedded in Euclidean space of $d$ dimension, (which is taken to be the
background spacetime). We therefore assume henceforth that the
world-sheets are closed and orientable surfaces, and attempt to work
out the contribution to the path integral for each fixed genus.

Fix for reference a closed oriented smooth surface
$\Sigma$ of genus $h$ (number of handles), and consider an
arbitrary smooth embeeding of $\Sigma$ in ${\bf R}^d$:
$$
X \equiv (X^1, \ldots, X^d) : \Sigma \longrightarrow {\bf R}^d
\eqno(1.3)
$$

\noindent
and simultaneously consider an arbitrary Riemannian metric $g$ on $\Sigma$:
$$
ds^2 = g_{ij} d\sigma^i d\sigma^j \quad, \quad i,j= 1,2 \eqno(1.4)
$$

\noindent
Here $(\sigma^1, \sigma^2)$ are local smooth coordinates on $\Sigma$.
The image $X(\Sigma)$ is to be considered as a typical
vacuum-to-vacuum world-sheet (path) and the random metric
$((g_{ij}))$ (having nothing to do with the metric induced on
$\Sigma$ via the embedding into Euclidean space) is an extra dynamical
variable which also has to be summed over in the Polyakov string
theory.

The fundamental definition is the {\bf Polyakov action} for that
world-sheet and that assigned metric:
$$
S(X,g) = \int\int_{\Sigma} [g^{ab} \frac{\partial X^\mu}{\partial
\sigma^a} \frac{\partial X^\mu}{\partial \sigma^b}] \sqrt{g}
d\sigma^1 d\sigma^2 \eqno(1.5)
$$

\noindent
Clearly, the part of the integrand (i.e. Lagrangian) in
square brackets is a real-valued function on $\Sigma$, and it is
being integrated with respect to the area-element,
$d_g(vol) = \sqrt{g} d\sigma^1 d\sigma^2$, induced on
$\Sigma$ by the metric $g$.

\noindent {\bf Notation :} Summation conventions over all repeated
indices are in use in (1.4), (1.5), and subsequently. In (1.5), indices
$a$ and $b$ are summed over $1 \le a , b \le 2$ , and $\mu$ is summed
over $1 \le \mu \le d$. Furthermore, $\sqrt{g}$ signifies the density
$\sqrt{det(g_{ij})}$, and $((g^{ab}))$ denotes the inverse matrix to
$((g_{ij}))$, as usual.

The basic problem therefore is to analyse the functional integral :
$$
Z = \int_{\{g_{ij}\}} \int_{\{X\}} e^{-S(X,g)} DX.Dg
\eqno(1.6)
$$
\noindent
That will represent the basic ``partition function'', or
vacuum-to-vacuum amplitude, for string propagations over world-sheets
with $h$ handles.

The remarkable thing is that, after taking care of certain symmetries in
the Polyakov action, the above integral does have a sensible reduction
to a {\it finite dimensional integration} over the moduli space ${\cal M}_h$ of
complex structures on a genus $h$ surface, provided the spacetime
dimension $d$ equals 26. Our first purpose, therefore, is to explain
concisely to a mathematical audience how the action prescription (1.5)
leads to a canonical measure -- called naturally the Polyakov measure --
on each moduli space ${\cal M}_h$.

\noindent
{\bf Remark on the classical theory for (1.5):} For a fixed
choice of embedding $X$, one may enquire as to what is the ``best''
(extremal) metric for the action (1.5). It is easily derived that the
action is extremised precisely for the metric $g$ on $\Sigma$ induced from
the (Euclidean) target space ${\bf R}^d$ via the embedding $X$. Namely,
$g_{ij} = \sum^d_{\mu =1} \frac{\partial X^\mu}{\partial \sigma^i}
\frac{\partial X^\mu}{\partial \sigma^j}$, is the ``classical'' metric.
So Polyakov action tells us to give this induced metric the highest weight
but average over {\it all} metrics in the path-integral (1.6). Note that
this fact, about the classical metric being the one induced from the
target {\bf R}$^d$, remains true even if we choose an arbitrary
Riemannian metric $((G_{\mu \nu}))$ in the background spacetime {\bf R}$^d$;
here, of course, we replace the Lagrangian integrand in (1.5) by the more
general: $[g_{ab} \frac{\partial X^\mu}{\partial \sigma^a}
\frac{\partial X^\nu}{\partial \sigma^b} G_{\mu \nu}]$.
As we will see, however, for the corresponding partition function (1.6),
the integral over the embedding variables (for fixed $((g_{ij}))$) is
{\it not} any more an (infinite dimensional) ``Gaussian''
if $((G_{\mu \nu}))$ is non-flat. For our purposes, therefore, we
restrict to the case $((G_{\mu \nu}))$ = Euclidean, just as in (1.5).

\vspace{.3cm}
\noindent {\bf I.C. Symmetries of the Polyakov action:}

To analyze (1.6), we first need to note that Polyakov's action has
certain {\it symmetries:}
$$
S(X^\mu + c^\mu, g) = S(X^\mu, g), ~ {\rm any}~ c^\mu \in {\bf R}^d.
\eqno(1.7)
$$
$$
S(f^\star X^\mu , f^\star g) = S(X^\mu, g), ~ {\rm any}~f \in Diff^+(\Sigma).
\eqno(1.8)
$$
$$
S(X^\mu, e^\phi g) = S(X^\mu , g), ~ {\rm any}~\phi \in
C^\infty(\Sigma,{\bf R}).
\eqno(1.9)
$$

\noindent
(1.7) corresponds to the fact that the action remains unchanged
if the embedded surface is simply translated in {\bf R}$^d$.
(1.8) says that, since (1.5) is invariantly defined, independent of
choice of coordinates on $\Sigma$, if we use any (orientation preserving)
diffeomorphism of $\Sigma$ to pullback both the metric and the embedding,
the action remains unperturbed. [Explicitly, $f^\star X^\mu = X^\mu \circ f$,
and $f^\star g$ is the metric on $\Sigma$ which assigns to any curve the
length that $g$ assigns to the $f$-image of the curve.] Finally, {\it
(1.9) is the truly non-trivial symmetry}, and says that, for a fixed
embedding, the value of {\it the action depends only on the
conformal class of the metric} $g$. [Conformal scaling,
$g \longmapsto (scaling function)g$ is called a
``Weyl-scaling'' by physicists.] Verification of (1.9) is immediate
since the Weyl factor cancels off between the square-bracketed
Lagrangian and the area-element term.

Clearly then, if the path integral (1.6) were actually computed over
{\it all} embeddings $X$ and {\it all} Riemannian metrics $g$, then one
would be getting infinite answers simply because one is {\it
overcounting} by (a) the ``volume of {\bf R}$^d$'', corresponding to the
arbitrary $c^\mu$ of (1.7) ; (b) the ``volume of Diff($\Sigma$)'' because
of (1.8) ; (c) the ``volume of positive functions (conformal factors) on
$\Sigma$'' because of (1.9). In other words, we can only expect to make
sense of (1.6) by quotienting out these symmetries -- namely by
integrating on the quotient space :
$$
\{ Embeddings\} \times \{Metrics\} /
\{{\bf R}^d \times Diff^+(\Sigma) \times Conf(\Sigma)\}
\eqno(1.10)
$$

\noindent {\bf Notation:} Write $Emb(\Sigma) = \{ X^\mu\}$ for the space
of all (smooth) embeddings of $\Sigma$ in {\bf R}$^d$, and $Met(\Sigma)$
for the space of all (smooth) Riemannian metrics on $\Sigma$.

\vspace{.3cm}
\noindent {\bf I.D. The integral over $\{X^\mu\}$} :

For every {\it fixed} $g
\in Met(\Sigma)$ on $\Sigma$ we will show below that, in analogy with
Gaussian (multivariate normal distribution) integrals, the integral over
\{ Embeddings\} / {\bf R}$^d$ can be carried out to produce a reasonable
answer. One should consider this section as providing
heuristic motivation for the following :

\noindent {\bf Proposition/Definition I.1.:}{\it For every fixed $g$ in
$Met(\Sigma)$, the inner integral in (1.6)is assigned the value:
$$
\int_{Emb(\Sigma)} e^{-S(X,g)} DX = \left[ \frac{det'
(-\Delta_g)}{\int_\Sigma d_g(vol)} \right]^{- d/2}
\eqno(1.11)
$$
\noindent
where $\Delta_g$ is the Laplace-Beltrami operator on functions
on $\Sigma$, and $det'(A)$ denotes the determinant of an
operator $A$ after discarding any zero eigenvalues.

The Polyakov integral (1.6) therefore becomes:
$$
Z = \int_{Met(\Sigma)}\left[\frac{det'(-\Delta_g)}{\int_\Sigma d_g(vol)}
\right]^{- d/2}Dg
\eqno(1.11*)
$$
with respect to a suitable volume element $[Dg]$ on the space of
Riemannian metrics on $\Sigma$.}

\noindent
{\it Note:} Determinants will be computed by heat-kernel (zeta-function)
regularization.

Since $(\Sigma,g)$ is a Riemannian manifold, we introduce the standard
$L^2$ inner-product on functions on $\Sigma$ by
$$
< f_1, f_2> = \int\int_{\Sigma} (f_1 f_2) d_g vol
\eqno(1.12)
$$

\noindent and obtain the Hilbert space $L^2(\Sigma)$ of square-integrable
real-valued functions with respect to this scalar product. The
Laplace-Beltrami operator, $\Delta_g$, is the formally self-adjoint
operator defined on sufficiently smooth functions on $\Sigma$ by :
$$
\Delta_g(f) = \frac{1}{\sqrt{g}} \left[ \frac{\partial}{\partial
\sigma^a} \left( \sqrt{g}^{ab} \frac{\partial f}{\partial \sigma^b}
\right) \right]  \eqno(1.13)
$$

\noindent
{\bf Lemma I.2 :} {\it The Polyakov action (1.5) can be rewritten as
$$
S(X,g) = \sum^d_{\mu =1} < X^\mu , (- \Delta_g) X^\mu >
\eqno(1.14)
$$
Notice that $(- \Delta_g)$ is a positive operator.}

\noindent {\bf Proof :} Triangulate $\Sigma$ so that each triangle falls
in a typical $(\sigma^1, \sigma^2)$ coordinate patch. Then integrating
by parts in (1.5) with respect to $\sigma^a$, on any triangle, gives an
integral over all boundary of the triangle plus the term
$$
- \int\int_{triangle} X^\mu \frac{1}{\sqrt{g}} \left[
\frac{\partial}{\partial \sigma^a} ({\sqrt{g}} g^{ab}
\frac{\partial X^\mu}{\partial \sigma^b}) \right]
\sqrt{g} d\sigma^1 d\sigma^2
$$

\noindent
Summing up over all the triangles, the boundary terms cancel
off and we are left with (1.14).
$\hfill{\Box}$

Now we are ready to give some heuristic arguments for (1.11).
Since all the $d$
target coordinates are on an equal footing, the equation (1.14) shows
that $\int e^{-S(X,g)} DX$ will be the product of $d$ identical
integrals :
$$
\int e^{< X^\mu, \Delta_g X^\mu >} DX^\mu ,~\mu=1,2,\ldots, d.
$$

Thus we are reduced to motivating the equation:
$$
\int_{\{Y:\Sigma \rightarrow {\bf R}\}} e^{- <Y, -\Delta_g Y>} DY = \left[
\frac{det' (-\Delta_g)}{Area(\Sigma,g)} \right]^{-1/2} \eqno(1.11')
$$
\noindent
This is clear. Indeed, suppose $\{e_o,e_1,e_2, \ldots\}$ is an
orthonormal basis for $L^2(\Sigma)$ consisting of eigenfunctions of
$-\Delta_g$. Since the $\Delta_g$-harmonic functions are the constants,
we have a 1-dimensional space of eigenfunctions with eigenvalue
$\lambda_o = 0$ generated by the constant function $e_o$. We set
$(-\Delta_g) e_k = \lambda_k e_k,$ with $\lambda_k > 0$ for all
$k=1,2,\ldots.$

Expanding in Fourier series an arbitrary $Y:\Sigma \rightarrow {\bf R}$ in
terms of the basis, we set
$$
Y = \sum^\infty_{k=0} y_k e_k
$$
\noindent with $ y_k = < Y, e_k >$.

Then one has
$$
-< Y, - \Delta_g Y > = - \sum^\infty_{k=0} |y_k|^2 {\lambda}_k
\eqno(1.15)
$$

If we take the functional measure ``$DY$'' in (1.11') to mean $dy_o dy_1
dy_2 dy_3 \ldots$ with Lebesgue measure in each factor, we see that we
have on our hands {\it a Gaussian integral}. Ignoring the $dy_o$
integral for the time being, this should produce $(\lambda_1 \lambda_2
\lambda_3 \ldots)^{-1/2}$ by comparison with the standard multivariate
normal integrals. Thus we make the reasonable definition that :
$$
\int e^{-\sum^{\infty}_{k=1} y^2_k \lambda_k} \left[ \Pi^\infty_{k=1} dy_k
\right] = \left[ det' (-\Delta_g) \right]^{- 1/2}
\eqno(1.16)
$$

Now, the integral over $y_o$, which should simply be ``the volume of
$y_o$-space'', clearly contributes an infinity which we wish to understand
and cancel against the infinity produced by the translation symmetry of
embeddings expressed in equation (1.7). We will explain this.

Since the eigenvector $e_o$ (corresponding to $\lambda_o = 0$) is a
constant, the normalization $\|e_o\|^2 = 1$ shows that the value of
$e_o$ is (Area $(\Sigma, g)$)$^{-1/2}$. In order to perform the
$y_o$-integral uniformly over all choices of the underlying metric $g$, we
stipulate the normalization $\int e^{-\| y_o \cdot 1\|^2} dy_o =1$.
[$\| \cdot \|$ always denotes the $L^2$ norm from (1.12). The
``$1$'' inside the norm means the (unnormalized) eigenvector given by the
constant (harmonic) function $1$. Note also that we ignore
the $\sqrt{2\pi}$ factors that appear in finite dimensional
Gaussian integrals.]

We were interested in analyzing
$\int e^{-\| y_o \cdot e_o \|^2 \lambda_o} dy_o$, where we think of
$\lambda_o$ as a {\it small} (but positive) eigenvalue
that we will ultimately send to zero. By the
normalized integral above, and the value of $e_o$, we are forced to the
conclusion that the integral above must behave like $(\lambda_o
e_o^2)^{-1/2} = \lambda_o^{-1/2} (Area (\Sigma,g))^{1/2}$.

So putting all the integrals over $(y_o, y_1, y_2 \dots)$ together
produces :
$$
\lambda_o^{-1/2} (Area (\Sigma, g))^{1/2} (det' (-\Delta_g))^{-1/2}
\eqno(1.17)
$$

Now remember that the symmetry exhibited in (1.7); in fact,
$Y \mapsto Y+c$, affects only the $y_o$-variable and produces the
overcounting infinity mantioned before. We think of this infinity as
``cancelling off'' the $\lambda_o^{-1/2}$ infinity appearing in (1.17)
(since $\lambda_{o} \ra 0$). Thus (1.17) becomes (1.11'),
which in turn motivates (1.11) itself.

In any event, we are only claiming that (1.11) is a well-motivated {\it
definition}, and a little thought about the above arguments
shows that it is the definition that is mathematically natural.

\vspace{.3cm}
\noindent {\bf I.E. Met$(\Sigma)$ and ${\cal M}_h$ :}

Taking care of the symmetries enjoyed by the action law, the Polyakov
integral (1.11*) has now metamorphosed to the following (still
rather enigmatic) shape :
$$
Z =
\frac{1}{vol(Diff^+(\Sigma)) \times vol(Conf(\Sigma))}
\int_{Met(\Sigma)} \left[ \frac{det'(-\Delta_g)}{Area(\Sigma, g)}
\right]^{- d/2}Dg
$$

$$
= \int_{Met(\Sigma) / Diff^+(\Sigma) \times Conf(\Sigma)}
{\left[\frac{det'(-\Delta_g)}{Area(\Sigma,g)} \right]^{-d/2}} Dg
\eqno(1.18)
$$

The final space over which one is now supposed to be integrating
is nothing other than the moduli space ${\cal M}_h$ parametrizing all
complex analytic structures on $\Sigma$. In fact, two metrics $g_1$ and
$g_2$ in $Met(\Sigma)$ are equivalent in ${\cal M}_h$ provided $g_2 =
f^\star (e^{\phi}g_1)$, for some $e^\phi \in Conf(\Sigma)$ and some
$f \in Diff^+ (\Sigma)$. Thus the complex structures, $\tau_1 $ and $\tau_2$,
obtained via isothermal parameters from the metrics $g_1$ and $g_2$,
respectively, are then biholomorphically equivalent via the biholomorphism
$f:\Sigma_{\tau_2} \longrightarrow \Sigma_{\tau_1}$. So we define:
$$
{\cal M}_h = Met(\Sigma) / \{Diff^+(\Sigma) \times Conf(\Sigma)\}
\eqno(1.19)
$$

The fundamental question, therefore, is whether the integrand in (1.18)
produces a well-defined measure on ${\cal M}_h$ for some choice of
the free parameter $d$ (the spacetime dimension).
We will sketch in the next sections how that pleasant state of affairs
transpires exactly when $d=26$.

\vspace{.3cm}
\noindent
{\it {\bf II:} POLYAKOV VOLUME FORM ON THE MODULI SPACE}

\noindent {\bf II.A. The Riemannian structure of $\Met$}

Consider the infinite dimensional manifold of Riemannian metrics,
$\Met$, on the fixed smooth surface $\S$. At the metric $g \in \Met$ we
can assign the following natural inner product to the tangent space of
$\Met$ thereat:
$$
<\delta g^{1}_{ab}, \delta g^{2}_{cd}> =
\int \int_{\S}(g^{ac}g^{bd} {\delta g^{1}_{ab}}{\delta g^{2}_{cd}})~
d_{g}(vol)
\eqno(2.1)
$$
where ${\delta g^{1}_{ab}}, {\delta g^{2}_{cd}}$ -- two symmetric
2nd-order covariant tensors on $\S$ --  are small perturbations of $g$,
representing two tangent vectors to $\Met$. Thus $\Met$ itself
qualifies as a Riemannian manifold (of infinite dimension).
{\it The Polyakov path integral, (1.11*) or (1.18), is to be interpreted
as an integration over $\Met$  with respect to the volume form induced
on $\Met$ by this Riemannian structure (2.1).}

\noindent
{\it Remark:} It is possible to make the mathematical considerations
regarding (2.1) completely rigorous by working with the Hilbert
manifold $\Met$ that one obtains by considering all metrics belonging
to a suitable Sobolev class $H^{s}$. Since the finite-dimensional
mathematics on the Teichm\"uller space that we finally arrive
upon is independent of these technical subtleties, we will not say
more about this matter in this exposition. In what follows it is possible
to restrict oneself to metrics, conformal scalings and diffeomorphisms
that are all of class $C^{\infty}$ on $\Sigma$.

As we know, the infinite dimensional manifold, $\Met$, is acted on by
the infinite dimensional group $\G = Diff^+(\Sigma) \times Conf(\Sigma)$,
(this is actually a semi-direct product), producing the quotient space
$\Mh$  --  an orbifold of finite dimension $(6h-6)$.

\noindent
{\it Remark:} If one assigns $\infty$ as the `dimension' of the space
of (local) functions on $\S$, then the choice of a Riemannian metric
tensor involves essentially three arbitrary functions -- so $\Met$ is
$3\infty$ dimensional. In that sense $Conf(\Sigma)$ is $\infty$
dimensional, and $Diff^{+}(\Sigma)$ is $2\infty$ dimensional. Thus,
the trading off of the parameters in $\Met$ for those in the gauge
group leaves a residual finite number of ``moduli parameters'' --
leading to the interesting equation ``$3\infty - 3\infty = 6h-6$''!

We shall work with the universal covering space of $\Mh$ in the
orbifold covering sense; that is the {\it Teichm\"uller space},
${\cal T}(\S) = \Th$. To define it, replace the group $Diff^+(\Sigma)$
in the above action by its identity component $Diff_{0}(\Sigma)$
(comprising diffeomorphisms homotopic to the identity). We obtain:
$$
{\cal T}_h = {Met(\Sigma)} /\{Diff_{0}(\Sigma) \times Conf(\Sigma)\}
\eqno(2.2)
$$
\noindent
Concomitantly, we shall denote the quotient projection from $\Met$
to $\Th$ by:
$$
\P:\Met \rightarrow \Th
\eqno(2.3)
$$

Notice that the quotient ${\Met}/Conf(\S)$ is precisely the space
of (smooth) Beltrami coefficients on $\S$ (see [N1]).

\noindent
{\bf Representing $\Th$ as a slice within $\Met$:}
We henceforth assume that the genus $h$ is at least two; (the
situation for spheres is trivial, and for tori the case is special and
easily treated.)  Then $\Th$ is a smooth manifold of real dimension
$(6h-6)$ (with a natural complex manifold structure); the discrete
{\it ``mapping class group'' ( or ``modular group'')},
$MCG(\Sigma) = \MCh = {Diff^{+}(\Sigma) / Diff_{0}(\Sigma)}$,
acts biholomorphically and proper discontinuously on $\Th$ producing
the quotient $\Mh$ as the (complex analytic) orbifold.

The space $\Th$ can be concretely pictured as the space of conjugacy
classes of ``marked'' Fuchsian groups $\{\Ga\}$ which are cocompact
and which produce quotient surfaces of the given genus. A Fuchsian
group $\Ga$ is marked by the choice of an isomorphism of the
fundamental group of $\S$ onto it. The moduli space $\Mh$ is just the
set of conjugacy classes of these Fuchsian groups (without markings).
See [N1] for this basic material.

Any smooth section (=right-inverse) of the quotient map $\P$ will be
called a {\it slice} in $\Met$.

\vspace{-.2cm}
\begin{center}
{\em Juliet:``What's in a name? that which we call a rose,}\\
{\em By any other name would smell as sweet;''}\\
\end{center}
\vspace{-.2cm}

Thus a slice is an embedded copy of Teichm\"uller space, $\Th$, in
$\Met$. Slices are $(6h-6)$ dimensional submanifolds of $\Met$, transverse
to the orbits of the ``gauge group'' $\G$. Mathematically speaking,
a slice represents the variation of the conformal moduli as variation of
Riemannian metric; in physics one says that choosing a slice ``fixes
the gauge freedom''.

\noindent
{\bf The Poincare slices and Weil-Petersson:}
The uniformization theorem guarantees that every Riemann surface
structure on $\S$ arises from a Poincare (hyperbolic) metric
of constant negative curvature ($-1$), that metric
being uniquely determined up to an arbitrary diffeomorphism. Define
therefore the following subset (infinite dimensional submanifold) of
$\Met$:
$$
Hyp(\S) = \{g \in \Met: curvature(g) \equiv -1\}
\eqno(2.4)
$$
The quotient ${\Met}/{Conf(\S)}$ is therefore in natural bijection
with the above submanifold of hyperbolic metrics.

Using the uniformization theorem with parameters, (see [N1]), we can
choose a smoothly varying family of hyperbolic metrics $\ga(t)$ in
$Hyp(\S)$ representing the Teichm\"uller space. Any such slice we
call a ``Poincare slice''. Thus $Hyp(\S)/{Diff_{0}(\S)}$ is the
Teichm\"uller space, and we can choose special Poincare slices that
are {\it orthogonal with respect to (2.1)} to the orbits of the gauge
group $Diff_{0}(\S)$. A convenient name we will adopt for such a slice is
{\it ``horizontal Poincare slice''}. The fundamental metric (2.1),
restricted to any horizontal Poincare slice, gives a Riemannian
structure to $\Th$ that is known classically as the {\it Weil-Petersson
metric} on the Teichm\"uller space. See Section III.A for more in this
direction.

\noindent
{\it Remark:} Having fixed a Poincare slice $\ga(t)$, any other
Poincare slice is then given by $f^{*}_{t}(\ga(t))$, where the
$f_t \in {Diff_{0}(\S)}$ are an arbitrary family of diffeomorphisms
that are chosen to depend smoothly on $t \in \Th$. It is therefore
clear that the metric (2.1) {\it cannot} induce the same metric on $\Th$
via an arbitrarily chosen Poincare slice; however, (2.1) does induce
Weil-Petersson on the {\it horizontal} Poincare slices defined above.

A natural way to select a Poincare slice is to utilize the unique
{\it harmonic diffeomorphism} (Eells-Sampson), that exists in the homotopy
class of the identity, between any two hyperbolic metrics on $\S$.
One may pullback the target metric via this diffeomorphism to obtain a
specific choice of hyperbolic metric on $\S$ representing the
Teichm\"uller class of the target metric. This Poincare slice is horizontal
(because harmonic Beltrami coefficients are orthogonal to the $Diff_{0}$
directions). See [J],[N3],[W] and the references therein.

\noindent
{Facts about (2.1):} {\it (a) The surface diffeomorphisms,
$Diff^{+}(\Sigma)$, act on $\Met$ as isometries of the
Riemannian structure (2.1).

\noindent
(b) However, the action of the Weyl-rescalings on $\Met$ do not
enjoy this compatibility with the metric (2.1).}

\noindent
{\it Remarks:} These are easily established.
Fact(b) above, namely that conformal rescalings fail to
preserve the metric on $\Met$, is at the root of the ``conformal
anomaly'' that we will have to grapple with in the material
below.

\vspace{.3cm}
\noindent {\bf II.B. Change of coordinates in $\Met$}

In order to analyse the integral (1.18) restricted to any slice, we must
clearly understand the nature of the metric (2.1) on $\Met$ in
coordinates that are along the gauge orbits, and complementary
coordinates in the Teichm\"uller directions along the slice.
Indeed, we want to factor out of the Polyakov integral the integrals
over the orbits of $Diff^+(\S)$ and $Conf(\S)$, as explained in (1.18).
Therefore it is natural to want to express the integral (1.18) as an
iterated integral over these gauge orbits and over the slice.

So fix once and for all some real analytic coordinates $((t))=((t_i))$,
$i$ running from $1$ to $N=6h-6$, on $\Th$; (the $((t_i))$
can be chosen as the real and imaginary parts of {\it holomorphic}
coordinates, and this can even be done globally over Teichm\"uller
space). We shall identify $\Th$ as this domain in $((t))$-space,
whenever convenient.

Let us {\it fix} a slice $K$, by choosing a right-inverse of $\P$:
$$
\ga:\Th \rightarrow \Met, ~~ \ga(\Th) = K
\eqno(2.5)
$$
and note that an arbitrary metric $\r \in \Met$ has a unique expression:
$$
\r = f^{*}[e^{\phi} {\ga}((t_i))]
\eqno(2.6)
$$
Here $((t_i))$ represents the Teichm\"uller point to which $\r$
projects by $\P$, $e^{\phi}$ is a conformal rescaling, and $f \in
Diff_{0}(\S)$. Thus the new coordinates for $\r$ are
$(f, \phi,((t_i))) \in Diff_{0}(\S) \times Conf(\S) \times \Th$.
Rescale the entire slice of metrics by the fixed $e^{\phi}$, and set:
$$
g = g((t)) = e^{\phi} {\ga}((t))
$$
$g((t))$ is an associated gauge-fixing slice, and note that
$f$ becomes an isometry from the metric $\r$ to the the metric $g$
because $\r=f^{*}(g)$. We remark that since $Diff_{0}(\S)$ is
the arc-component of the identity in $Diff^+(\S)$, it is possible
to express $f$ as $exp(\xi)$ where $\xi$ is a smooth vector-field on $\S$.

We need to compute the Riemannian structure (2.1) of $\Met$ in these new
coordinates. At any given point $\r \in \Met$, we must understand small
metrical variations $\d \r_{ab}$ in terms of changes in these new
$Gauge \times \Th$ coordinates. Instead of working with a small change
of $\r$ we will work with a corresponding small change of $g$;
remembering (Fact (a) above), that the pullback action by the fixed
diffeomorphism $f$ is an isometric automorphism of $\Met$, we lose
nothing by this.

To this end we first write the arbitrary small change in
the metric in the form:
$$
\d\r = f^{*}(\d g)
\eqno(2.7)
$$
utilizing the {\it same} diffeomorphism $f$ as in (2.6).
The definition of $\d g$ implies that it involves: a diffeomorphism
close to the identity (which we write as $exp(\xi)$), a small
change in the Weyl rescaling, $\d\phi$, and small changes $((\d t_i))$
in the Teichm\"uller coordinates;
$\d g = (exp \xi)^{*}[e^{\phi + \d\phi}{\ga}((t+\d t))] - g$

We will work in the tangent space to $\Met$ at the point $g$. Recall
that the {\it trace} (with respect to $g_{ab}$) of a symmetric
tensor $\s_{ab}$ is by definition the contraction $g^{ab}\s_{ab}$.
Observe that the traceless symmetric second-order covariant tensors,
$dh_{ij}$, and the pure-trace tensors of the form $\d \phi g_{ab}$,
constitute {\it orthogonal} spaces with respect to the fundamental
inner product (2.1). It is therefore convenient to express the
general metrical deformation $\d g$ as a sum of traceless and pure-trace
parts. That will be done below.

\noindent
{\it Notation:} We shall sum over repeated indices in the formulae of
this article. The roman letters $a, b$ etc will usually vary over the
surface coordinates (i.e., $1 \le a,b \le 2$), whereas the indices
$i$, $j$, $m$, $n$ will usually run through Teichm\"uller coordinates and
thus range from $1$ to $N=(6h-6)$. Indices will be raised and lowered
with respect to $g$. $D(\xi)$ will denote covariant derivative of the
vector field $\xi$, with respect to $g$.

\noindent
{\bf Lemma II.1:} {\it The infinitesimal change $\d g$ defined above
decomposes as:
$$
\d g_{ab} = \d\phi g_{ab} + (P\xi)_{ab} + \d t^{i} T^{i}_{ab}
\eqno(2.8)
$$
with the last two terms constituting the tracefree part.
Here
$$
P:\{Vector~fields~on~\S\} \rightarrow
\{Symmetric~traceless~ (2,0)~tensors~on~\S\}
$$
is the 1st order elliptic operator (depending on $g$), given by
$$
(P\xi)_{ab} = D_{a}(\xi_b) + D_{b}(\xi_a) - g_{ab}D_{c}(\xi^c)
$$
and,
$$
T^{i}_{ab} = {\frac {\partial}{\partial t_i}}[g_{ab}((t))] -
{\frac{1}{2}}g^{cd}{\frac{\partial}{\partial t_i}}[g_{cd}((t))]
$$
The $T^i((t))$ are the (trace free parts of) the tangent vectors
in the Teichm\"uller directions to the gauge slice $g((t))$. }

\noindent
{\bf Proof:} Differential geometry on the surface gives:
$(exp \xi)^{*}g_{ab} = g_{ab} + D_{a}(\xi_b) + D_{b}(\xi_a)$;
and by Taylor expansion one sees:
$ e^{\d\phi}g_{ab} = g_{ab} + \d\phi g_{ab} + o(\d \phi)$.
Computing with the help of these gives (2.8). In (2.8), by suitably
redefining $\d\phi$, all the pure-trace component of the $\d g$
has been absorbed in the first term on the right hand side.
$\hfill{\Box}$

Now, the various spaces of tensor-fields on the surface $(\S,g)$ carry
natural inner products induced by the metric $g$. For the space of
vector-fields the pairing is explicitly:
$$
(\xi,\eta)_{g} = \int\int_{\S}[\xi^{a}\eta^{b}g_{ab}]d_{g}vol
\eqno(2.9a)
$$
and for the covariant 2nd-order traceless tensors the pairing is:
$$
(R,S)_{g} = \int\int_{\S}[R_{ab}S_{cd}g^{ac}g^{bd}]d_{g}vol
\eqno(2.9b)
$$
(Whenever the metric $g$ is clear from the context we will take the
liberty of suppressing that subscript.)

We can think of the operator $P$ as an unbounded closed operator
defined on the appropriate dense domain of the Hilbert space of vector
fields. Then, by basic functional analysis, the target Hilbert space of
symmetric 2nd order tracefree tensors decomposes into the orthogonal
direct sum:
$\{2nd~order~symmetric~tracefree~tensors\}=\{Range~P\}\oplus\{Ker~P^{*}\}$.

The range of $P$ constitutes the piece of the tangent space to $\Met$
arising from pullbacks of the metric $g$ by diffeomorphisms close to
the identity. The Weyl (conformal) scalings of $g$ are already
absorbed in the pure-trace part of the infinitesimal deformation
$\d g$. Consequently, the residual piece given by $KerP^{*}$
is the finite-dimensional part comprising tangent vectors in the
Teichm\"uller (slice) directions.

The final upshot is that the tangent space at $g$ to the infinite
dimensional space $\Met$ decomposes as an orthogonal direct sum:
$$
T_{g}(\Met) = \{pure~trace\}\oplus \{Range~P\}\oplus \{Ker~P^{*}\}
\eqno(2.10)
$$

In order to understand the Teichm\"uller deformations piece, which is
of central interest to us, we need therefore to analyse the adjoint
of $P$ (computed with respect to the inner products (2.9a) and
(2.9b)). One easily verifies that $P^{*}$ is nothing other than the d-bar
operator on the space of (smooth) quadratic differentials on $(\S,g)$
(thought of as a Riemann surface). {\it Consequently,
the kernel, $Ker P^{*}$, is the $6h-6$ dimensional (real) subspace
of the tracefree 2nd order symmetric tensors consisting of the
(real and imaginary parts of) holomorphic quadratic differentials.}

\noindent
{\bf Teichm\"uller's Lemma: $T^{*}_{t}(\Th) \equiv
H^{0}(X_t,K(X_{t})^{2})$:}
The upshot of the discussion above will be recognized by those
familiar with Teichm\"uller theory as a variant of ``Teichm\"uller's
lemma''; that lemma describes the complex cotangent space (at the
conformal structure $g$) to the Teichm\"uller space $\Th$,
as precisely (namely, canonically isomorphic to) the vector space
of holomorphic quadratic differentials on $(\S,g)$. See [N1]. Indeed,
the tangent space at $t \in \Th$ to Teichm\"uller space, is, by
Kodaira-Spencer, canonically isomorphic to the first cohomology
$H^{1}(X_t,K^{-1})$, which, by Serrre duality, is canonically
isomorphic to the dual space of global holomorphic sections of $K^2$.
By the Riemann-Roch theorem, this last space is a complex vector space
of complex dimension $(3h-3)$ (for each $t \in \Th$).
(Note: $K=K(X_{t})$ denotes the holomorphic cotangent bundle of the
Riemann surface $X_t$.)

\noindent
{\bf Bases for $H^{0}(X_t,K^{2})$ as $X_t$ varies in $\Th$:}
To facilitate our work,
we therefore introduce an arbitrary auxiliary choice of basis for the
holomorphic quadratic differentials on each of the Riemann surfaces
$X_t = (\S, g((t)))$ -- i.e., all over the Teichm\"uller space $\Th$.
Thus let $Q^m((t))$, $1 \le m \le (6h-6)$, be a basis of the $6h-6$
dimensional (over reals) vector space of holomorphic quadratic
differentials on the Riemann surface $(\S,g((t)))$, where the $Q^m$
are assumed to vary smoothly with $((t)) \in \Th$. (That is easily
accomplished by, for instance, utilizing Poincare theta-series to
define the $Q^m$ with respect to real-analytically varying Fuchsian
groups $\Gamma_{((t))}$.) Notice the important fact that once we have
chosen bases as above on each surface along a particular slice, the job
is accomplished for {\it all} slices, -- because {\it the $Q^m$ depend
only on the conformal class} of the metric under consideration.

{\it Clearly, Teichm\"uller's lemma, asserts that each $Q^m((t))$ can be
considered as a nowhere vanishing $1$-form over $\Th$}. We will need
this interpreation in Lemma II.5 below.

Substituting (2.8) into the Riemannian structure formula (2.1) we
obtain immediately the desired re-expression for the Riemannian norm
on the tangent space to $\Met$ at $g$:
$$
||\d g||^{2} = ||\d\phi||^{2} + (P\xi, P\xi) +
(T^i,Q^m)((Q^m,Q^n))^{-1}(Q^n,T^j)\d t_i \d t_j
\eqno(2.11)
$$
We explain the notations:
The $L^2$ norm for functions on the surface $(\S,g)$ was already shown
in (1.12), -- and that defines the first term $||\d\phi||^{2}$.
The $T^i((t))$ are, as shown above, the (trace free parts of)
the tangent vectors in the Teichm\"uller directions.
Each $T^i$ is, of course, a second-order symmetric traceless
covariant 2nd-order tensor on $(\S,g)$. All the pairings appearing in the
second and third terms on the right of (2.11) are the inner products
exhibited in (2.9b). Thus, setting $R=S=P\xi$ in (2.9b) gives the
second term, and similarly for the pieces in the third term.
(Note that the inversion involved in the last term is matrix inversion.)

\noindent
{\bf Proposition II.2:} {\it The volume measure induced by the Riemannian
structure (2.1) on the infinite dimensional manifold $\Met$ expressed
in the coordinates $(\phi, \xi, ((t_i)))$ is:
$$
[Dg] = [(detP^{*}P)^{1/2}][det((Q^m,Q^n))]^{-1/2}[det(T^i,Q^m)]
[D\phi][D \xi]dt_1\wedge \cdots \wedge dt_{6h-6}
\eqno(2.12)
$$
By $[D\phi]$ we mean the volume measure on $Conf(\S)$
induced by the $L^2$-norm (1.12) on the space of rescaling functions
$\phi$. Similarly, $[D\xi]$ is the volume element on $Diff_{0}(\S)$
arising from the inner product (2.9a) on vector fields.}

\noindent
{\bf Proof:} The expression of the Riemannian structure (2.1) in the
form (2.11) is adapted precisely to the slice $K$. Thus the usual
formula $d_{g}vol = \sqrt{det((g_{ij}))}dx_1\wedge\cdots\wedge dx_{M}$
for the Riemannian volume element on a M-dimensional manifold, when
applied formally to this infinite dimensional $\Met$, gives (2.12).
Indeed, in the infinite dimensional pieces corresponding to the first
two terms on the right of (2.11), one needs to compute the
determinants of the operators $I=Identity$ and $P$, respectively. For
later convenience we have substituted $|det P| = [(detP^{*}P)^{1/2}]$.
$\hfill{\Box}$

Recall that the integrations over the gauge variables, $(\phi,\xi)$,
is expected to produce the volumes of the gauge groups that we have been
desiring to factor out (see (1.18)).  Therefore, ignoring
the $[D\phi][D \xi]$ in the above volume element of $\Met$ is
clearly a reasonable method of accomplishing that aim, -- and it is
the physicist's way of reducing the Polyakov integral to a finite
dimensional integral over the given slice. Dropping $[D\phi][D \xi]$
from (2.12) we obtain therefore {\it a volume element on the slice $K$}:
$$
d\mu_{K} =
{\left[\frac{(detP_{g}^{*}P_{g})}{det((Q^m,Q^n)_{g})}\right]^{1/2}}
[det((T^i,Q^m)_{g})] dt_1\wedge \cdots \wedge dt_{6h-6}
\eqno(2.13)
$$
We emphasize that as the metric $g=g((t))$ varies over $K$, the operator
$P=P_g$ varies, and thus the infinite determinant $[(detP^{*}P)^{1/2}]$
above is a non-trivial function of $((t))$.

Now we have reached the goal of expressing the Polyakov integral as a
well-defined finite dimensional integration over any chosen slice:

\noindent
{\bf II.3 The Polyakov prescription:} {\it In the light of Proposition
II.2 and the discussion above, we interpret the Polyakov integral (1.18)
as the following integral over the (arbitrarily chosen) slice:
$$
\int_{Slice~K}{\left[\frac{det'(-\Delta_g)}{Area(\Sigma,g)}
\right]^{-d/2}}d\mu_{K}
\eqno(2.14)
$$}
The integrand above is (for any choice of $d$) a well-defined volume
form on the slice $K$, and since $K$ is an embedded copy of $\Th$ the
volume form can be considered as living on $\Th$.

\noindent
{\it Remark:} We have been working over the Teichm\"uller space, $\Th$,
rather than over the moduli space, $\Mh$ -- as was prescribed by
(1.18). So, strictly speaking, the Polyakov integral (2.14) should be
divided by the ``cardinality of the mapping class group $\MCh$''.
Alternatively, we should think of the integral in (2.14) only over
a fundamental domain for the action of $MCG$. In any event, we are only
interested in the {\it volume form = integrand} on $\Th$ or $\Mh$,
rather than in the integral itself. The basic fact is that the Polyakov
measure on $\Th$ is $MCG$-invariant, but its total integral over $\Mh$
(i.e., over any fundamental domain) is infinite (see (3.3) below).

\vspace{.3cm}
\noindent {\bf II.C. Finally! The Polyakov measure on $\Th$:}

Now that we have reduced ourselves from the original Polyakov path
integral (1.6) to the integral (2.14), the fundamental question is
whether (2.14) is {\it independent of the choice of the slice}. The
main result is:

\noindent
{\bf Theorem II.4:} {\it When $d/2=13$, the integrand in (2.14), as
a volume form on the Teichm\"uller space $\Th$, is {\bf independent}
of the choice of the slice $K$. This is the Polyakov volume form,
$d(Poly)$,  on $\Th$.  It is mapping class group invariant, is absolutely
continuous with respect to the Lebesgue measure class on $\Th$.}

\noindent
{\bf Addendum to Theorem:} {\it d(Poly) has the following expression in
terms of the Weil-Petersson volume form on $\Th$:
$$
d(Poly) = [(Z'_{G}(1))^{-13}Z_{G}(2)] d(Weil-Pet)
\eqno(2.15)
$$
Here $Z_{G}(s)$ is the Selberg zeta function (an entire function)
associated to the (variable) Fuchsian group $G$. }

\noindent
{\it Remark:} Throughout the above considerations, we are ignoring
overall constant factors in expressions for the Polyakov form.
It is our attitude that the Polyakov volume form
on each genus moduli is an interesting mathematical object defined
only up to an arbitrary scaling constant. See (3.3) below also.

The notations in (2.15) will be explained as we go along.
Combining (2.13) and (2.14) we realize that we have to deal with the
following volume form on the slice $K$:
$$
dP_{K} =
{\left[\frac{det'(-\Delta_g)}{Area(\Sigma,g)}\right]^{-d/2}}
{\left[\frac{(detP_{g}^{*}P_{g})} {det((Q^m,Q^n)_{g})}\right]^{1/2}}
[det((T^i,Q^m)_{g})] dt_1\wedge \cdots \wedge dt_{6h-6}
\eqno(2.16)
$$
(the metric $g$ varies over $K$.)

The aim is to understand how far the above volume density depends on
the choice of the Teichm\"uller slice. Let us break up the work into
natural pieces:

\noindent
{\bf Lemma II.5:} {\it The volume form on $\Th$ given by:
$$
d\lambda = [det((T^i,Q^m)_{g})] dt_1\wedge \cdots \wedge dt_{6h-6}
\eqno(2.17)
$$
is independent of the choice of the slice. (It depends only on the
choice of the bases of quadratic differentials $Q^m((t))$, which we
have fixed, once and for all, over the entire Teichm\"uller space $\Th$.)

Indeed, recalling that the $Q^m((t))$ can be considered as $1$-forms
over $\Th$, this volume form is the following top dimensional form:
$$
d\lambda = Q^1 \wedge \cdots \wedge Q^{6h-6}
\eqno(2.18)
$$}

\noindent
{\bf Proof:}
A new slice $K^\sharp$ is obtained from the slice $K$ by scaling the
metrics comprising $K$ in a $((t))$-dependent fashion:
$g^{\sharp}((t)) = e^{\phi((t))}g((t))$. The trace-free parts of the
tangent vectors to the slices are related by:
$T^{i\sharp} = e^{\phi}T^{i}$. Computing the relevant pairings now, by
the law (2.9b) (with respect to the metrics $g$ and $g^\sharp$,
respectively), we obtain
$$
(T^{i\sharp}, Q^m)_{g^{\sharp}} = (T^{i}, Q^m)_{g}
\eqno(2.19)
$$
Note that the above equality simply represents the fact that
the $T^i$ are tangent vectors to the Teichm\"uller space,
whereas the $Q^m$ are co-tangent vectors to the same space,
and that the pairing is just the {\it duality pairing} of
tangents with cotangents. Hence, it is not surprising that
(2.19) holds independent of scalings in the metric.

Consequently, (2.17) is independent of the choice of the
slice, as claimed. The rest is straightforward.
$\hfill\Box$

Consider therefore the fundamental function $F$ defined below (for
each value of the unspecified ``space-time'' dimension $d$). Each
$F=F_d$ is a real-valued function on the entire space of metrics, $\Met$:
$$
F_{d}(g) = {\left[\frac{det'(-\Delta_g)}{Area(\Sigma,g)}\right]^{-d/2}}
{\left[\frac{(detP_{g}^{*}P_{g})} {det((Q^m,Q^n)_{g})}\right]^{1/2}}
\eqno(2.20)
$$
We shall see that $F_{26}$ is singled out!:

\noindent
{\bf Proposition II.6:} {\it The function $F_d$ is invariant along the
orbits of the gauge group $\{Diff^{+}(\S) \times Conf(\S)\}$ if and
only if $d=26$; $F_{26}$ thus descends as a function to the
moduli space $\Mh$. Namely one has the remarkable invariance:
$$
F_{26}(f^{*}(e^{\phi}g)) = F_{26}(g)
\eqno(2.21)
$$
for all $g \in \Met$, under all conformal rescalings and pullbacks by
diffeomorphisms.}

\noindent
{\bf Proof:}
Let us first note the significance of the two infinite dimensional
determinant factors that constitute the function $F_d$. Define:

$$
D_0(g)  = {\left[\frac{det'(-\Delta_g)}{Area(\Sigma,g)}\right]}
\eqno(2.22)
$$

$$
D_1(g) = {\left[\frac{(detP_{g}^{*}P_{g})}{det((Q^m,Q^n)_{g})}\right]}
\eqno(2.23)
$$
Remember that the operator $\Delta_g$ is the Laplacian
acting on smooth functions on $\S$, and that $P_{g}^{*}P_{g}$
is the Laplacian acting on smooth vector-fields.

The problem is to understand the variation of these two determinantal
functions on $\Met$ under a conformal rescaling of the metric. This is
a standard type of problem for the application of heat kernel
techniques -- the results are clearly exposed in Alvarez[Alv]. See
also [AN]. Recall the heat-kernel (i.e., zeta-function) regularized
determinant for operators in infinite dimensions:
$$
log(det' D) = - {\lim_{\epsilon \rightarrow 0}}
\int_{\epsilon}^{\infty}[Trace'(e^{-tD})]{\frac{dt}{t}}
\eqno(2.24)
$$
The notation $det'$ and $Trace'$ mean that the determinant
or trace is being calculated after discarding the zero eigenvalues --
namely on the orthogonal complement of the kernel of $D$.

Consider an arbitrary infinitesimal conformal scaling of the metric
$g$ to $\ga = e^{\d\phi}g$. The corresponding perturbation in the
non-zero eigenvalues of the relevant operators appearing in (2.22) and
(2.23) can be calculated. The first variations turn out
as below, ($R_{g}$ denotes the scalar curvature for metric $g$):

$$
\d(log D_{0}) = {\frac{1}{12\pi}}\int_{\S}R_{g}\d\phi d_{g}(vol)
\eqno(2.25)
$$

$$
\d(log D_{1}) = {\frac{13}{12\pi}}\int_{\S}R_{g}\d\phi d_{g}(vol)
\eqno(2.26)
$$

Therefore, $\d F_{d}$ is ${\frac{1}{12\pi}}(13-{\frac{d}{2}})$ times
the integral appearing on the right of the formulae (2.25) and (2.26).
The invariance property (2.21) for $F_{26}$ follows immediately.
$\hfill\Box$

\noindent
{\bf Proof of Theorem II.4:}
Combining the gauge-invariant function $F_{26}$ on $\Met$ with the
slice-independent volume form $d\lambda$ of Lemma II.5, we see that we
have proved that {\it the Polyakov volume form on $\Th$ has the
expression}:
$$
d(Poly) = F_{26}d\lambda
\eqno(2.27)
$$
Since the bases for quadratic differentials were introduced simply as
an artifice to help simplify our formulae, it is obvious that the
volume form $d(Poly)$ above is also independent of the choice of
these bases $\{Q^m((t))\}$. [Note: Glancing at (2.16), we may see directly
that the combination of $Q$-dependent pieces:
$[det((T^i,Q^m)_{g})]{[det((Q^m,Q^n)_{g})]^{-1/2}}$, evidently
remains invariant when the $Q^m$ are subjected to an arbitrary
(in general $t$-dependent) change of basis.]

We have therefore established the existence of the Polyakov volume form
on $\Th$, independent of all arbitrary choices, as stated in Theorem II.4.
$\hfill\Box$

As general references for the matter we presented see [Alv],[AN],[Nel].
The expression for $d(Poly)$ in terms of $d(Weil-Pet)$ will be
dealt with in the next section.



\vspace{.3cm}
\noindent
{\it {\bf III:} GEOMETRY OF TEICHM\"ULLER SPACE AND POLYAKOV VOLUME}

\noindent
{\bf III.A. Weil-Petersson and Polyakov:}

\vskip .3cm
\noindent
{\bf Weil-Petersson Riemannian structure on $\Th$:} Since the
cotangent space to $\Th$ at any point $t$ is canonically isomorphic to
the vector space of holomorphic quadratic differentials on the Riemann
surface $X_{t} = (\S,\ga(t))$, a hermitian structure is induced on
$\Th$ by utilizing the Petersson pairing of holomorphic quadratic
differentials on $X_t$. In fact, let $\ga(t)$ be a hyperbolic metric
representing the conformal structure $t \in \Th$; if the metric is
expressed as $\ga = \lambda(z)|dz|$ in terms of local holomorphic
(=isothermal) coordinates, then set
$$
(Q^1,Q^2)_{WP} = \int_{\S}{Q^1(z)}{\overline{Q^2(z)}}{\lambda(z)}^{-2}
dz \wedge d\bar{z}
\eqno(3.1)
$$
This is Weil-Petersson hermitian metric of $\Th$. See [IT],[J],[N1]. The
corresponding volume form on $\Th$ we shall denote by $d(Weil-Pet)$.

\noindent
{\bf Lemma III.1:} {\it (a) The Riemannian structure (2.1) of $\Met$,
restricted to a horizontal Poincare slice, is Weil-Petersson on $\Th$.

\noindent
(b) Choose the bases $\{Q^m((t))\}$ for the spaces of holomorphic
quadratic differentials to be orthonormal (pairing (2.9b)) with
respect to the corresponding hyperbolic metrics $\ga((t))$.
Then the (slice-independent) volume form
$d\lambda = Q^1 \wedge \cdots \wedge Q^{6h-6}$ described in
Lemma II.5 is $d(Weil-Pet)$.}

\noindent
{\bf Proof} (a):
A holomorphic quadratic differential $Q$ on $(\S,\ga)$ is a cotangent
vector to $\Th$ there; the Weil-Petersson pairing sets up an
isomorphism of the cotangent space with the tangent space -- thus $Q$
corresponds to a certain tangent vector $v_Q$ to $\Th$. Identifying
$\Th$ with the slice, the tangent $v_Q$ is given by an infinitesimal
change in the hyperbolic metric $\ga$. But that perturbation of metric,
say $\d\ga_{Q}$, is, in its turn, represented by a certain Beltrami
differential on $(\S,\ga)$. Indeed, one sees that the co-vector $Q$
corresponds to the harmonic (Bers') Beltrami differential
$\mu_{Q} = \overline{Q(z)}{y^2}$, ($z=x+iy$).

Above, we are utilizing uniformization to represent the
Riemann surface $(\S,\ga)$ as $H/{\Ga}$, $H$ being the Poincare
upper half-plane, and $\Ga$ a Fuchsian group operating thereon.
The hyperbolic metric $\ga = {y^{-1}}|dz|$ therefore gets deformed to
$\ga + \d\ga = {y^{-1}}|dz + \varepsilon\mu(z)d\bar{z}|$ with
$\mu = {\mu}_{Q}$ as shown. See [N1] for details.

Expanding out the right side we may compute the $\d\ga = {\d\ga}_{Q}$;
then substituting in the formula (2.1) we get
$({\d\ga}_{Q^1}, {\d\ga}_{Q^2}) = (Q^1,Q^2)_{WP}$, as desired.

\noindent
(b): Clearly, at every point $((t)) \in \Th$, the co-vectors
$\{Q^m((t))\}$ now constitute an orthonormal basis for
$T^{*}_{((t))}(\Th)$ with respect to the Weil-Petersson inner-product.
Therefore, $Q^1 \wedge \cdots \wedge Q^{6h-6} = d(Weil-Pet)$
$\hfill\Box$

\noindent
{\it Remark:} The bases of quadratic differentials above can be
constructed using Poincare theta series with respect to varying
Fuchsian group $\Ga((t))$ and then applying the Gram-Schmidt
orthonormalization procedure. Since the bases are arbitrary up to any
smooth choice $\Th \rightarrow GL(6h-6,{\bf R})$, we note that the
corresponding measure $Q^1 \wedge \cdots \wedge Q^{6h-6} = d\lambda$
can be, a priori, chosen quite arbitrarily.

\vskip .3cm
\noindent
{\bf Selberg zeta and determinants for Laplacians:}  Associated to the
hyperbolic Riemann surface $(\S,\ga) = H/{\Ga}$, $\Ga$ the
uniformizing (purely hyperbolic, cocompact) Fuchsian group, Selberg
defines
$$
Z_{\Ga}(s) = {\prod}_{\a}{\prod}_{n \geq 0}[1-exp\{-l(g_{\a})(s+n)\}]
\eqno(3.2)
$$
where $\{g_{\a}\}$ are a complete set of representatives for the
conjugacy classes of the primitive elements of $\Ga$ (i.e., those
allowing no non-trivial $n$-th roots in $\Ga$); $l(g)$ denotes the length
of the geodesic in the free homotopy class determined by $g$. This
function is defined by the above product for $Re(s)>1$, and can be
shown to have an analytic continuation as an entire function in the
whole $s$-plane.

The Selberg trace formula, that goes hand in hand with the
above zeta function, allows one to express the heat-kernel regularized
determinants of the Laplace operators on $(\S,\ga)$  -- on functions,
vector-fields and higher order tensor fields -- as certain values of
the above holomorphic function. As explained, for example, in
[DP], one gets $det'(-\Delta)=Z'(1)$ (up to a constant independent of
the Fuchsian group), and similarly, $det(P^{*}P)=Z(2)$. (All the
numbers $Z(k)$, integral $k$, have interpretations as determinants of
Laplacians acting on suitable spaces of tensors on $(\S,\ga)$.)

\noindent
{\bf Proof of Addendum to Theorem II.4:}
Glancing now at (2.22), (2.23), and armed with (b) of Lemma III.1, we
have clearly obtained the expression for $d(Poly)$ claimed in
the Addendum to Theorem II.4. (Note: Over hyperbolic metrics the
Gauss-Bonnet shows that $Area(\S,\ga)=4\pi(h-1)$; moreover, the determinant
appearing in the denominator of (2.23) we had already normalized to
unity by our choice of orthonormal bases for quadratic differentials.)

\noindent
{\it Remark:}
The asymptotics of the Selberg zeta functions, when the Fuchsian groups
approach the boundary (Deligne-Mumford boundary) of $\Mh$, can be
used to show that the Polyakov volume form blows up near the boundary
so fast that the fundamental amplitude integral is divergent:
$$
\int_{\Mh}d(Poly) = \infty
\eqno(3.3)
$$
(Recall that the Weil-Petersson volume of $\Mh$, on the other hand, is
finite.) (3.3) demonstrates that no finite answer is forthcoming via
the Polyakov prescription for the ``vacuum-to-vacuum'' amplitude of
the bosonic string.

\vspace{.3cm}
\noindent
{\bf III.B. Complex geometry and Polyakov:}

We are now at the point where we can explain our deepest reason for
excitement about the Polyakov measure. {\it In fact, $d(Poly)$ has a
natural and simple construction (over each $\Th$) by simply the complex
analytic geometry of $\Th$ -- with no inputs whatsoever from physics.}

\noindent
{\bf The universal family and allied bundles over $\Th$:}

The Teichm\"uller spaces are fine moduli spaces. Namely,
the total space $\S \times {\cal T}_h$ admits a natural complex
structure such that the projection to the second factor
$$
{\psi}_h:{\cal C}_h:=\S \times{\cal T}_h \rightarrow {\cal T}_h
\eqno{(3.4)}
$$
gives the universal Riemann surface over ${\cal T}_h$.
This means that for any $\eta \in {\cal T}_h$, the
submanifold $\S \times\eta$ is a complex
submanifold of ${\cal C}_h$, and the complex structure on $\S$ induced
by this embedding is represented by $\eta$. As is well-known, (Chapter 5
in [N1]), the family $\Ch \rightarrow \Th$ is the {\it universal} object
in the category of holomorphic families of genus $h$ marked
Riemann surfaces.

We recall now the {\it determinant of cohomology}
construction for obtaining line bundles over any parameter space of
Riemann surfaces. Let $X$ be a compact Riemann surface and $L$ be a
holomorphic line bundle on $X$.  The determinant of cohomology of $L$ is
defined to be the $1$-dimensional complex vector space $detH^{*}(L) =
det(L)$:
$$
det(L) = (\ext{top}H^0(X,L))\bigotimes (\ext{top}H^1(X,L)^*)
\eqno(3.5)
$$

We are ready to pass to families. Let $\pi:\X \longrightarrow S$
be a holomorphic family of compact Riemann
surfaces parametrized by a base $S$ which is complex analytic space.
The prototypical instance is the
universal family described over $\Th$. The  definition demands that
$\pi$ be a holomorphic submersion with compact and connected fibers
of complex dimension one -- each fiber being of genus $h$.

For any point $s \in S$, the construction (3.5) gives a
complex line $det(L_s)$. The basic fact is that these lines fit
together to give a holomorphic line bundle on $S$ (see [BGS]), which
is called the {\it determinant of cohomology bundle} of $L_S$, and is
denoted by $det(L_S)$. (Note: This is not to be confused with the top
exterior power of $L_S$ -- which is a line bundle over the total space
$\X$.) The general definition arises from the direct image functors
$R_{*}$ of algebraic geometry. See [KM], [D].

Over each genus Teichm\"uller space we thus have a sequence of natural
determinant of cohomology bundles arising from the powers of the relative
tangent bundles along the fibers of the universal curve. Indeed,
let ${\om}_h \longrightarrow {\cal C}_h$ be the relative cotangent bundle
for the projection ${\p}_h$ in $(3.4)$. The determinant line bundle
over $\Th$ arising from its $n$-th tensor power is fundamental for us,
and we shall denote it by:
$$
DET_{n,h}:= det({\om_{h}}^{n}) \longrightarrow \Th, ~~n \in \Z
\eqno(3.6)
$$
Applying Serre duality shows that there is a canonical isomorphism
$DET_{n,h} = DET_{1-n,h}$, for all $n$. $DET_{0,h} = DET_{1,h}$ is
called the {\it Hodge} line bundle over $\Tg$.

\noindent
{\it Remarks:}
(i) The determinant of cohomology bundle, $det(L_S)$, is functorial
with respect to base change. Given any morphism
$\ga: S' \longrightarrow S$ this means, in particular,
that there is a canonical isomorphism of the determinant
of cohomology associated to the pulled back family over $S'$ with the
pullback by $\ga$ of the original determinant bundle over $S$. See
[BNS] for details.

\noindent
(ii) The determinant of cohomology construction, $det(L_S)$, produces a
bundle over the parameter space $S$ induced by the bundle over
the total space $\X$; now, the Grothendieck-Riemann-Roch (GRR)
theorem (see [BGS], [D]) gives a canonical isomorphism of $det(L_S)$
with a combination of certain bundles (on $S$) obtained from the direct
images of the bundle $L_S$ and the relative tangent bundle
$T_{\X/S}$. The GRR theorem is important for our work below.

\noindent
(iii) Whenever the bundles $L_S$ and the relative tangent bundle above
are assigned smooth hermitian structures, the determinant of cohomology
inherits a smooth hermitian structure, due to Quillen [Q], in a
functorial way. See [BGS], [D], [BNS]. The canonical isomorphisms
arising from base change (see (i)) then become {\it isometric
isomorphisms}.

The vertical (=relative) tangent bundle for the universal family,
namely the bundle $T_{{\Ch}/{\Th}} \rightarrow \cal C$ comprising the
tangent bundles of the fibering Riemann surfaces, carries
a canonical hermitian structure by assigning the Poincar\'e
hyperbolic metric on each surface (utilizing the uniformization theorem
with moduli parameters). As a consequence, applying remark (iii) above,
the holomorphic line bundles $DET_{n,h}$ carry natural {\it Quillen
hermitian metric} arising from the Poincar\'e metrics on
the fibers of $\Ch$.

Observe that by the naturality of the above constructions it follows
that the action of $\MCh$ on $\Th$ has a natural lifting as automorphisms
of these $DET$ bundles. These automorphisms are unitary with respect to
the Quillen structure.

\noindent
{\bf The fibers of $DET_{n,h}$:} Notice that the fiber of Hodge,
over $X \in \Th$, is the top exterior power of
the vector space of holomorphic Abelian differentials on $X$. Similarly,
the fiber of $DET_{2,h}$ is the $(3h-3)$ exterior power of the vector
space of holomorphic quadratic differentials on $X$, and so on. By some
complex geometry (e.g., Poincare theta series with respect to
holomorphically varying quasi-Fuchsian groups), one can create the
natural {\it $MCG$ equivariant holomorphic vector bundle}, $V_{n,h}$
over $\Th$ by attaching over $X$ the fiber $H^{0}(X,K_{X}^{n})$
($n=1,2,..$). Therefore, (3.6) shows that the top exterior
powers of these vector bundles over $\Th$ are nothing other than the
determinant of cohomology bundles $DET_{n,h}$ we are describing.

In particular, by Teichm\"uller's lemma we know $V_{2,h}$ is
canonically isomorphic to the holomorphic cotangent bundle of
Teichm\"uller space. Thus, there are canonical isomorphisms:
$$
\ext{(3h-3)}{V_{2,h}} = \ext{top}{T^{*}{\Th}} = DET_{2,h}
\eqno(3.7)
$$
each of the above being a description of the canonical line bundle
over $\Th$.

We also recall that holomorphic $1$-forms can be paired naturally, in
the $L^2$ sense, on any Riemann surface $X$ via the ``Hodge pairing''
$$
(\a,\b)=-i\int\int_{X}{\a \wedge \bar{\b}}
\eqno(3.8)
$$
This gives a $MCG$ invariant hermitian metric on $V_{1,h}$, and hence
induces a $MCG$ invariant hermitian structure (called the Hodge metric)
on the Hodge bundle. (The Hodge metric has to be modified suitably by
a factor given by the determinant of the Laplacian to obtain the Quillen
metric on that bundle.)

\noindent
{\bf Mumford isomorphisms and Polyakov volume:}
These natural line bundles over $\Th$ will be considered as
$MCG$-equivariant line bundles, and the isomorphisms we talk about will
be $MCG$-equivariant isomorphisms. By applying the
Grothendieck-Riemann-Roch theorem indicated above, Mumford [Mum] had
shown that $DET_{n,h}$ is canonically isomorphic to a certain fixed
({\it genus-independent!}) tensor power of the Hodge bundle.

\noindent
{\bf Proposition III.2:}
{\it
$$
DET_{n,h} \cong {DET_{1,h}}^{\otimes(6n^2-6n+1)}
\eqno{(3.9)}
$$
The above isomorphism is isometric with respect to the Quillen metrics.
An isomorphism (3.9) is ambiguous only up to a non-zero scalar.}

\noindent
{\it Remarks:} (i) It can be shown that the Picard group of (isomorphism
classes of) all $\MCh$-equivariant holomorphic line bundles over $\Th$,
(which can be identified as the group of holomorphic line bundles over the
moduli space $\Mh$), is a cyclic group generated by the Hodge bundle.
This can be shown from work of Harer, Mumford and Arbarello-Cornalba.
That cyclic group is of order $10$ for $h=2$, and is infinite cyclic
for $h > 2$. The above Proposition shows that the natural
sequence of determinant of cohomology line bundles pick out a certain
expicit subsequence from the Picard group.

\noindent
(ii) The $MCG (\equiv \MCh)$ invariant holomorphic
functions on $\Th$ -- namely those that descend to $\Mh$ --
are {it constant}, for all $h>2$. This follows from
the Satake compactification theory for $\Mh$. Thus, although $\Mh$ is
non-compact, in certain aspects it behaves like a compact analytic space.
Consequently, (at least when $h > 2$), any two isomorphisms between bundles
over $\Mh$ can only be ambiguous up to a scale factor, as asserted.

There is a remarkable connection, discovered by Belavin and Knizhnik [BK],
between the Mumford isomorphism above for the case $n=2$, [i.e., that
$DET_{2}$ is the $13$-th tensor power of Hodge], and the existence of the
Polyakov string measure on the moduli space $\Mg$. In fact this $13$ is
the same lucky number as the value of $d/2$ in Proposition II.6. See
[BK],[Bos],[Nel].

First, it is elementary to see that assigning a hermitian metric
on the canonical bundle of any complex space gives a measure on that
space. Indeed, fixing a volume density on a space simply amounts to
fixing a fiber metric on the canonical line bundle, $K$, -- because then
we know the unit vectors in $K$ -- and absolute square gives volume.
Now, $K$ for the Teichm\"uller space is nothing other than $DET_{2,h}$
-- recall (3.7) above. Thus, any natural hermitian structure on $DET_{2,h}$
becomes a choice of a natural volume form on $\Th$.

\noindent
{\bf Theorem III.3:}
{\it The Hodge bundle, $Hodge \equiv DET_{1,h}$, has its natural Hodge
metric, arising from (3.8). We may transport the corresponding metric on
${Hodge}^{13}$ to $DET_{2,h}$ by Mumford's isomorphism, (the choice
of isomorphism being unique up to scalar) -- thereby obtaining a $MCG$
invariant volume form, say $d(Mum)$, on $\Th$ (also unambiguous up to
the choice of a scalar). This $d(Mum)$ is none other than the
Polyakov volume form $d(Poly)$ of Theorem II.4.}

\noindent
{\bf Sketch of Proof:}
Choose a local holomorphic frame, $(\om_{1}(t),\cdots,\om_{h}(t))$ for the
rank $h$ vector bundle $V_{1,h}$ over a local holomorphic
$t$-patch in $\Th$. Suppose that $(\p_{1}(t),\cdots,\p_{3h-3}(t))$
is a local holomorphic frame for $V_{2,h}$ over the same
neighbourhood such that under the Mumford isomorphism:
$$
[\om_1 \wedge\cdots\wedge\om_{h}]^{13} \mapsto
(\p_1\wedge\cdots\wedge \p_{3h-3})
\eqno(3.10)
$$
Then the Mumford volume form on the $t$-coordinate patch is
immediately seen to have the expression:
$$
d(Mum) = [det((\om_{i}(t),\om_{j}(t)))]^{-13}{d\lambda}
\eqno(3.11)
$$
with the matrix of pairings being the Hodge pairings of (3.8), and
where the measure $d\lambda$ is given by Lemma II.5, equation (2.18):
$d\lambda = \p_1\wedge\cdots\wedge \p_{3h-3}\wedge
\overline{\p_1}\wedge\cdots\wedge \overline{\p_{3h-3}}$

Utilizing therefore the bases $Q^m$ for quadratic differentials given
by precisely these $\p_k$ and their conjugates, we may compare
$d(Mum)$ of (3.11) with $d(Poly)$ as shewn in equation (2.20); we
derive that the {\it ratio} of these two volume forms on $\Th$ is the
following function:
$$
\frac{d(Poly)}{d(Mum)} =  F_{26}(g(t))[det((\om_{i}(t),\om_{j}(t)))]^{13}
$$
$$
= {\left[\frac{det'(-\Delta_g)}{Area(\Sigma,g)}\right]^{-13}}
{\left[\frac{(detP_{g}^{*}P_{g})} {det((\p_m,\p_n)_{g})}\right]^{1/2}}
[det((\om_{i}(t),\om_{j}(t)))]^{13}
\eqno(3.12)
$$
Above, $g=g(t)$ is any Riemannian metric representing the conformal
structure $t \in \Th$.

Since both volume forms under consideration are, by construction,
modular invariant, we note that the above ratio, say $G(t)$, is a real
valued $MCG$ invariant function on the Teichm\"uller space.

\noindent
{\it Claim:}
$\frac{d(Poly)}{d(Mum)} = G(t)$ is the absolute value square of
a global {\it holomorphic} function $f$ on $\Th$.

Basically one computes $\d\bar{\d}[log G(t)]$ and shows that this is
identically zero. This local computation is explained in
[BK],[Nel] and other references. Thus $ holomorphic f = f_U$ exists
on any local $(t)$-neighbourhood $U$ satisfying $G(t) = |f(t)|^2$.
But $\Th$ is simply connected, and it is clear that these
local $f_U$ can be patched up to analytically continue along all paths
in $\Th$. It follows directly that there is a holomorphic function $f$
defined on the entire $\Th$ such that $G(t) \equiv |f(t)|^{2}$.

Finally, we show that this global $f$ is necessarily modular invariant.
In fact, under any modular transformation $f$ can only be multiplied by
factor of absolute value one  -- thus $f$ gives rise to a character
for the modular group. But, since that group, $MCG$, is generated by
elements of finite order (!), which all have fixed points in their
action on $\Th$, it follows that the character must be trivial. That
is as desired.

But such $f$, and hence $G$, must be constant (recall the Satake
compactification remark) -- and we are through.
$\hfill\Box$

\noindent
{\it Remarks:} (i) The expression for $d(Mum)$ in (3.11) also allows one to
study the blow up of the volume form as one approaches the
Deligne-Mumford boundary of the moduli space. This can be used to reprove
the divergence of the string amplitude stated in (3.3).

\noindent
(ii) As we said above Theorem III.3, any natural hermitian structure
on $DET_{2,h}$ is the assignment of a natural volume form on $\Th$.
Therefore the Quillen hermitian structure on this $DET$ bundle
assigns also a volume element to Teichm\"uller/moduli
space. This is being studied in relation to the
Polyakov volume. For example, is the Quillen volume of $\Mh$ finite?

\vspace{.3cm}
\noindent
{\it {\bf IV:} THE UNIVERSAL TEICHM\"ULLER SPACE OF COMPACT SURFACES}

The moral of the story above is that the presence of the {\it Mumford
isomorphisms} over the moduli space of genus $g$ Riemann surfaces
describes the Polyakov measure structure thereon! It is consequently a
natural problem to try to create {\it genus-independent} versions
of these constructions by working over some universal parameter space
that parametrizes Riemann surfaces of varying genus. Indeed, the
string may sweep over world-sheets of arbitrary genus -- so that the
problem raised is fundamental to non-perturbative string theory.

In this and the following section we give a concise report of our
work [BNS] where we succeed in finding a genus-independent description
of the Mumford isomorphisms over a certain universal parameter space $\Tin$.
We would like to emphasize one point: all our constructions are
equivariant under the action on $\Tin$ of an intriguing new mapping
class group $MC_{\infty}$.

\vspace{.3cm}
\noindent
{\bf IV.A. The direct limit, $\Tin$, of classical Teichm\"uller spaces:}

We start with a fundamental topological situation. Let
$$
\pi:\tX \longrightarrow X
\eqno(4.1)
$$
be an {\it unramified covering map}, orientation preserving,
between two compact connected oriented two manifolds $\tX$ and $X$ of
genera $\tg$ and $g$, respectively. Assume $g \geq 2$. The degree of the
covering $\pi$, which will play an important role, is the ratio
of the respective Euler characteristics; namely, $deg(\pi)=(\tg-1)/(g-1)$.

Given any complex structure on $X$, we may pull back this structure
via $\pi$ to a complex structure on $\tX$.
Now, the homotopy lifting property guarantees that there is a unique
diffeomorphism ${\tilde f}\in{\rm Diff}_0({\tX})$ which is
a lift of any given $f \in{\rm Diff}_0({X})$.
Mapping $f$ to $\tilde f$ defines an injective
homomorphism of ${\rm Diff}^{+}_0(X)$ into ${\rm Diff}^{+}_0({\tX})$.
Consequently, $\pi$ induces an injection of the smaller Teichm\"uller
space into the larger one:
$$
{\cal T}(\pi):{\cal T}_g \longrightarrow {\cal T}_{\tg}
\eqno(4.2)
$$
It is known that this map ${\cal T}(\pi )$ is a {\it proper holomorphic
embedding} between these finite dimensional complex manifolds; $\Tpi$
respects the quasiconformal-distortion (=Teichm\"uller) metrics.
{}From the definitions it is evident that this embedding between the
Teichm\"uller spaces depends only on the (unbased) isotopy class of
the covering $\pi$.

At the level of Fuchsian groups, one should note that any covering
space $\pi$ corresponds to the choice of a subgroup $H$ of finite index
(=$deg(\pi)$) in the uniformizing group $G$ for $X$, and the embedding
(4.2) is then the standard inclusion mapping $\T(G) \rightarrow \T(H)$;
(see Chapter 2, [N1]).

\noindent
{\it Remark:} One notices that $\cal T$ is a
contravariant functor from the category of closed oriented
topological surfaces, morphisms being covering maps, to the
category of finite dimensional complex manifolds and holomorphic
embeddings. We shall have more to say along these lines below.

We construct a category $\A$ of certain topological objects and
morphisms: the objects, $Ob(\A)$, are a set of compact oriented topological
surfaces each equipped with a base point ($\star$), there being exactly
one surface of each genus $g \geq 0$; let the object of genus $g$ be
denoted by $X_g$. The morphisms are based isotopy classes of pointed
covering mappings
$$
\pi: (X_\tg, \star) \rightarrow (X_g, \star)
$$
there being one arrow for each such isotopy class. {\it Note that the
monomorphism of fundamental groups induced by (any representative of)
the based isotopy class $\pi$, is unambiguously defined.}

\noindent
{\bf The direct system of classical Teichm\"uller spaces:}
Fix a genus $g$ and let $X = X_g$.
Observe that all the morphisms with the fixed target $X_g$:
$$
M_g = \{\a \in Mor(\A): Range(\a)=X_g \}
\eqno(4.3)
$$
constitute a {\it directed set} under the partial ordering given by
factorisation of covering maps.
Thus if $\a$ and $\b$ are two morphisms from
the above set, then $\b \succ \a$ if and only if the image of the
monomorphism $\pi_1(\b)$ is contained within the image of $\pi_1(\a)$;
that happens if and only if there is a commuting triangle of morphisms:
$\b = \a \circ \theta$. It is important to note that the factoring
morphism $\theta$ is {\it uniquely determined} because we are working
with base points.
[Remark: Notice that the object of genus $1$ in $\A$ only has
morphisms to itself -- so that this object together with all its
morphisms (to and from) form a subcategory.]

By (4.2), each morphism of $\A$ induces a proper, holomorphic,
Teichm\"uller-metric preserving embedding between the corresponding
finite-dimensional Teichm\"uller spaces. We can thus create the natural
{\it direct system of Teichm\"uller spaces} over the above directed
set $M_g$, by associating to each $\a \in M_g$ the Teichm\"uller space
$\T(X_{g(\a)})$, where $X_{g(\a)} \in Ob(\A)$ denotes the domain surface
for the covering $\a$. To each $\b \succ \a$ one associates the
corresponding holomorphic embedding $\T(\theta)$ (with
$\theta$ as in the paragraph above). From this direct system
we form the {\it direct limit Teichm\"uller space over $X=X_g$}:
$$
\Tin(X_g) = \TinX := ind. lim. \T(X_{g(\a)})
\eqno(4.4)
$$
$\TinX$ is a metric space with the Teichm\"uller metric, and it also
has a natural Weil-Petersson Riemannian structure obtained from scaling
the Weil-Petersson pairing on each finite dimensional stratum, $\Tg$,
by the factor $(g-1)^{-1}$. See [NS].

$\Tin$ is our {\it commensurability Teichm\"uller space} -- which
will serve as the base space for universal Mumford isomorphisms.

\noindent
{\bf The Teichm\"uller space, $\THin$, of the hyperbolic solenoid:}
Over the same directed set $M_g$ we may also define a
natural {\it inverse system of surfaces}, by associating to $\a \in
M_g$ a certain copy, $S_{\a}$ of the pointed surface $X_{g(\a)}$.
[Note: Fix a universal covering over $X=X_g$. $S_{\a}$ can be taken to
be the universal covering quotiented by the action of the subgroup
$Im(\pi_1(\a)) \subset {\pi_1}(X,\star)$.] If $g \ge 2$, then
the inverse limit of this system is the {\it universal solenoidal
surface} $H_{\infty}(X) = proj.lim.X_{g(\a)}$ that was studied in [S],[NS].

In fact, $\Hin$ is a compact topological space that fibers over $X$
with the fibers being Cantor sets. The path components of $\Hin$, (with
leaf-topology), are simply connected two-manifolds restricted to each of
which the projection $\pi_{\infty}:\Hin \rightarrow X$ becomes a
universal covering. There are uncountably many of these path components
(``leaves'') in $\Hin$, and each is a dense subset in $\Hin$.
Each leaf is thus, morally speaking, a hyperbolic plane. That is why
we call $\Hin$ the universal hyperbolic solenoid.
The facts above follow from a careful study of this inverse system of
surfaces, the main tool being the lifting of paths in $X$ to its coverings.

As explained in [S],[NS], the solenoid $\Hin$ has a natural
Teichm\"uller space comprising equivalence classes of complex structures
on the leaves -- the leaf complex structures being required to vary
continuously in the fiber (Cantor) directions. In particular, any
complex structure assigned to any of the surfaces $X_{g(\a)}$
appearing in the inverse tower can be pulled back to all the surfaces
above it -- and therefore assigns a complex structure of the sort
demanded on $\Hin$ itself. These complex structures that arise from
some finite stage can be characterized as the ``transversely locally
constant'' (TLC) ones (see [NS]), and they comprise precisely the
{\it dense} subset $\TinX$ sitting within the separable Banach manifold
$\THinX$. We collect some of these thoughts in the:

\noindent
{\bf Proposition IV.1:} {\it The ``ind-space'' $\TinX$, (see [Sha]), is
the inductive limit of finite dimensional complex manifolds, and hence
carries a complex structure defined strata-wise. The completion of
$\TinX$ with respect to the Teichm\"uller metric is the separable
complex Banach manifold $\THinX$.

Alternatively, $\TinX$ can be embedded in Bers' universal Teichm\"uller
space, $\T(\Delta)$, as a directed union of the Teichm\"uller spaces of
Fuchsian groups. (The Fuchsian groups vary over the finite index
subgroups of a fixed cocompact Fuchsian group $G$, $X=\Delta/G$.) The
closure in $\T(\Delta)$ of $\TinX$ is a Bers-embedded copy of $\THinX$.}

\noindent
{\it Remark:}
It is evident, but important to note, that the spaces $\TinX$ and
$\THinX$ we are dealing with do {\it not} really depend on the choice of $X$.
If we were to start with a surface $X'$ of different genus (also greater
than one), then we could pass to a common covering surface (always
available!), and hence the limit spaces we construct would be isomorphic.

\vspace{.3cm}
\noindent
{\bf IV.B. The commensurability mapping class group $\MCin = Aut(\Tin)$:}

A remarkable but obvious fact about the above construction is that
every morphism $\pi:Y \rightarrow X$ of $\A$ induces a natural
{\it Teichm\"uller metric preserving homeomorphism}
$$
\Tin(\pi): \TinY \longrightarrow \TinX
\eqno(4.5)
$$
$\Tin(\pi)$ is invertible simply because
the morphisms of $\A$ with target $Y$ are cofinal
with those having target $X$ (thus all finite ambiguities are
forgotten in passing to the inductive limits!). It is also clear
that $\Tin(\pi)$ is a biholomorphic identification (with respect to
the strata-wise complex structures). [Note that $\Tin$ is covariant
-- whereas the Teichm\"uller functor $\T$ itself was contravariant.]

It follows that each $\TinX$, and so also its metric completion $\THinX$,
is equipped with a large {\it automorphism group} -- one from each undirected
cycle of morphisms of $\A$ starting from $X$ and returning to $X$.
By repeatedly using pull-back diagrams (i.e., by choosing appropriate
connected component of the fiber product of covering maps), it is fairly
easy to see that the automorphism of $\TinX$ arising from
any (many arrows) cycle can be obtained simply from a two-arrow cycle
$\tX {\ra \atop \ra} X$.

Namely, whenever we have (the isotopy class of) a ``self-correspondence''
of $X$ given by any two non-isotopic coverings, say $\a$ and $\b$,
$$
\tX {\ra \atop \ra} X
\eqno(4.6)
$$
we can create a corresponding automorphism $R \in Aut(\TinX)$
defined as the composition: $R = {\Tin(\b)}\circ{(\Tin(\a))^{-1}}$.

These automorphisms constitute a group that we shall call
the {\it commensurability modular group}:
$$
\MCin(X) = Aut(\TinX)
\eqno(4.7)
$$
acting on $\TinX$ and on $\THinX$.

To clarify matters further, consider the abstract graph
($1$-complex), $\Ga(\A)$, obtained from the topological
category $\A$ by looking at the objects
as vertices and the (undirected) arrows as edges. It is
clear from the definition above that the fundamental group
of this graph, viz. $\pi_{1}(\Ga(\A),X)$, is acting on $\TinX$
as these automorphisms. (We may fill in triangular $2$-cells
in this abstract graph whenever two morphisms (edges) compose to
give a third edge; the thereby-reduced fundamental group of this
$2$-complex also produces on $\TinX$ the action of $\MCin(X)$.)

\noindent
{\bf Making explicit the genus one situation:}
For the genus one object $X_1$ in $\A$, we know
that the Teichm\"uller space for any unramified covering
is but a copy of the upper half-plane $H$. The maps $\Tpi$ are
M\"obius identifications of copies of the half-plane with itself,
and we easily see that the pair $(\Tin(X_1),CM_{\ify}(X_1))$ is
identifiable as $(H,PGL(2,\QQ))$. Notice that the action has dense
orbits in this case. Anticipating for a moment the definition of $Vaut$
given below, we remark that $GL(2,\QQ) \cong Vaut(\Z \oplus \Z)$, and
$Vaut^{+}$ is precisely the subgroup of index $2$ therein, as expected.

On the other hand, if $X \in Ob(\A)$ is of any genus $g \geq 2$,
then we get an infinite dimensional ``ind-space'' as $\TinX$ with
the action of $\G(X)$ on it as described. Since the tower of coverings
over $X$ and $Y$ (both of genus higher than $1$) eventually become
cofinal, it is clear that {\it for any choice of genus higher than one
we get {\bf one isomorphism class} of pairs} $(\Tin, \MCin)$.

\noindent
{\bf Virtual automorphism group of $\pi_{1}(X)$ and $\MCin$:}
In the classical situation, the action of the mapping class group
$MCG(X)$ on $\T(X)$ was induced by the action of (isotopy classes of)
self-homeomorphisms of $X$; in the direct limit set up we  now have the
more general (isotopy classes of) self-correspondences of $X$ inducing
the new mapping class automorphisms on $\TinX$. In fact, we will see that
our group $\MCin$ corresponds to ``virtual automorphisms'' of the
fundamental group $\pi_{1}(X)$, -- generalizing exactly the classical
situation where the usual $Aut(\pi_{1}(X))$ appears as the action
via modular automorphisms on $\T(X)$.

Given any group $G$, one may look at its ``partial'' or ``virtual''
automorphisms; as opposed to usual automorphisms that are defined on all
of $G$, for virtual automorphisms we demand only that they be defined
on some finite index subgroup. To be precise,
consider isomorphisms $\r:H \ra K$ where $H$ and $K$ are
subgroups of finite index in $G$. Two such isomorphisms (say $\r_1$
and $\r_2$) are considered equivalent if there is a finite index subgroup
(sitting in the intersection of the two domain groups) on which they
coincide. The equivalence class $[\r]$  -- which is like the {\it
germ} of the isomorphism $\r$ -- is called a {\it virtual automorphism}
of $G$; clearly the virtual automorphisms of $G$ constitute a group,
christened $Vaut(G)$, under the obvious law of composition, (i.e., compose
after passing to deeper finite index subgroups, if necessary).

Clearly $Vaut(G)$ is trivial unless $G$ is infinite (though there do
exist infinite groups -- see [MT] -- such that $Vaut$ is trivial).  Also
evident is the fact that $Vaut(group)=Vaut(any~finite~index~subgroup)$.
Since we shall apply this concept to the fundamental group of a surface of
genus $g$, ($g>1$), the last remark shows that our $Vaut(\pi_{1}(X_g))$
{\it is genus independent}!

In fact, $Vaut$ presents us a neat way of formalizing the ``two-arrow
cycles'' (4.6) which we introduced to represent elements of $\MCin$.
Letting $G = \pi_{1}(X)$, (recall that $X$ is already equipped with a
base point), we see that the diagram (4.6) corresponds exactly to the
following virtual automorphism of $G$:
$$
[\r] = [{\b}_{*}\circ{\a}_{*}^{-1}:{\a}_{*}(\pi_{1}(\tX)) \rightarrow
{\b}_{*}(\pi_{1}(\tX))]
\eqno(4.8)
$$
Here ${\a}_{*}$ denotes the monomorphism of the fundamental group
$\pi_{1}(\tX)$ into $\pi_{1}(X) = G$, and similarly ${\b}_{*}$ etc..
We let $Vaut^{+}({\pi}_{1}(X))$ denote the subgroup of $Vaut$ arising
from pairs of orientation preserving coverings.

As we said, all the automorphisms arising from arbitrarily
complicated cycles of coverings (i.e., any finite sequence
of morphisms starting and ending at $X$), are each obtained
from these simple two-arrow pictures (4.6).
The reduction of any many-arrow cycle in $\Ga(\A)$ to a two-arrow cycle
utilizes successive fiber product diagrams; there is some amount
of choice in this reduction process, and one may obtain different
two-arrow cycles starting from the same cycle -- but the virtual
automorphism that is defined is unambiguous. The final upshot is:

\noindent
{\bf Proposition IV.2:} {\it One has natural surjective group
homomorphisms:
$$
{\pi}_{1}(\Ga(\A),X) \rightarrow Vaut^{+}({\pi}_{1}(X)) \rightarrow
Aut(\TinX) \equiv \MCin(X)
$$}

\noindent
{\it Remarks:} (i) The concept of $Vaut$ has already arisen in group
theory papers -- for example [Ma],[MT]. I am grateful to Chris Odden for
pointing out these references to me, and for interesting discussions.

(ii) There is a natural representation of $Vaut(\pi_{1}(X))$ in the
homeomorphism group of the unit circle $S^1$, by the standard theory of
{\it boundary homeomorphisms} (see, for example, [N1]).
This leads to many obvious questions whose
answers are not obvious. In recent work, Robert Penner and the present
author have proved an isomorphism between the group $\MCin$ above and a
direct limit of Penner's Ptolemy groups that operate on his ``Tessellations''
version of universal Teichm\"uller space. That isomorphism is connected
with the following conjecture, and we hope to report on these matters in
forthcoming publications.

\noindent
{\bf Topological transitivity of $\MCin$ on $\Tin$ and allied issues:}
Does $\MCin$ act with dense orbits in $\Tin$? That is a basic query. This
question is directly seen to be equivalent to the following old conjecture
which, we understand, is due to L.Ehrenpreis and C.L.Siegel:

\noindent
{\bf Conjecture IV.3:} {\it Given any two compact Riemann surfaces, $X_1$
(of genus $g_1 \geq 2$) and $X_2$ (of genus $g_2 \geq 2$), and given any
$\epsilon > 0$, can one find finite unbranched coverings $\pi_1$ and
$\pi_2$ (respectively) of the two surfaces such that
the corresponding covering Riemann surfaces $\tX_1$ and $\tX_2$ are of the
same genus and there exists a $(1+\epsilon)$ quasiconformal homeomorphism
between them. (Namely, $\tX_1$ and $\tX_2$ come $\epsilon$-close in the
Teichm\"uller metric.)}

\noindent
{\it Remark:} Since the uniformization theorem guarantees that the {\it
universal} coverings of $X_1$ and $X_2$ are {\it exactly} conformally
equivalent, the conjecture asks whether we can obtain high {\it finite}
coverings that are {\it approximately} conformally equivalent.

\vspace{.3cm}
\noindent
{\it {\bf V:} UNIVERSAL POLYAKOV-MUMFORD ON $\Tin$:}

The goal is to construct natural geometrical
fiber bundles over the commensurability Teichm\"uller space,
which when restricted to the finite dimensional
strata, $\Tg$, become the Hodge (or higher $n$) $DET$ bundles thereon.
The bundles over $\Tin$ should be related by universal Mumford
isomorphisms, which restrict to the strata as the finite dimensional
Mumford isomorphisms explained in (3.9). We are able to carry out this
entire program in a $\MCin$ equivariant fashion ([BNS]).

\vspace{.3cm}
\noindent
{\bf V.A. The Main Lemma:}
We invoke into play the arbitrary unramified finite covering
$\pi: \tX \rightarrow X$, and recall that $\Tpi$ (equation(4.2)) is
the associated holomorphic embedding of Teichm\"uller spaces.

The fundamental question is whether there exists any natural
relationship between the line bundle $\Dng$ over the Teichm\"uller
space $\Tg=\T(X)$ and the bundle ${\Tpi}^{*}{\Dnt}$ obtained as the
pullback of the corresponding determinant of cohomology bundle
over the larger Teichm\"uller space $\Tt=\T(\tX)$. For example, we are
asking, is there any natural relationship between the two Hodge
bundles?

We have an elegant answer to this question that forms the
foundation for our genus-independent description of Mumford
isomorphisms. In effect, $\Dng$ raised to the tensor power $deg(\pi)$,
simply extends naturally over the larger Teichm\"uller space $\Tt$ as
the $\Dnt$ bundle thereon! We prove this by utilizing the
Grothendieck-Riemann-Roch (GRR) theorem (of [D]) in crucial ways.

\noindent
{\bf Lemma  V.1:} {\it The two holomorphic line bundles (with Quillen
hermitian structures), $(\Dng)^{deg(\pi)}$ and ${\Tpi}^{*}{\Dnt}$, on
${\cal T}_g$ are canonically isometrically isomorphic for every integer $n$.
(The isomorphism is canonical up to the choice of a $12$th root of unity.)
In other words, there is a canonical isometrical line bundle morphism
$\Ga(\pi)$ lifting $\Tpi$ and making the following diagram commute:
$$
\matrix{
{\Dng}^{deg(\pi)}
&\mapright{{\Ga}(\pi)}
&{\Dnt}
\cr
\mapdown{}
&
&\mapdown{}
\cr
\Tg
&\mapright{\Tpi}
&\Tt
\cr}
$$
The maps $\Ga(\pi)$ are functorial, so that for any commuting triangle of
morphisms in the category $\A$, the corresponding $\Ga$-lifts also
commute.}

\noindent
{\bf Curvature forms of $DET$ bundles and Lemma V.1:}
The existence of the canonical relating morphism between the above
determinant bundles (fixed $n$) in the fixed covering
situation was first conjectured and deduced by us (see [BN]) utilizing the
differential geometry of the Quillen metrics. Recall that the
Teichm\"uller spaces $\Tg$ and $\Tt$ carry natural symplectic
forms -- the Weil-Petersson K\"ahler forms -- which are in fact
the {\it curvature forms of the Quillen metrics} of these $DET$
bundles ([Wol], [ZT], [BGS]). If the covering $\pi$ is unbranched of
degree $d$, a direct calculation shows that this natural WP
form  on $\Tt$ (appropriately renormalized by the degree $d$)
pulls back to the WP form of $\Tg$ by $\Tpi$. Equality of the curvature
forms leads one to expect the isomorphism of  ${DET_{n,g}}^{d}$
with ${\Tpi}^{*}{\Dnt}$.  This intuition is what is behind the more
sophisticated GRR proof of the Lemma above.

\vspace {.3cm}
\noindent
{\bf V.B. Power-law principal bundle morphisms over Teichm\"uller spaces}

We desire to obtain certain canonical geometric objects over the
inductive limit of the finite dimensional Teichm\"uller spaces
by coherently fitting together the determinant line bundles
$DET_{n,g}$ thereon. To this end it is necessary to find a canonical
mapping relating $\Dng$ itself to $\Dnt$ utilizing the Lemma above.

Now, given any complex line bundle $\la \ra T$ over any base $T$,
there is a certain {\it canonical} mapping of $\la$ to any positive
integral ($d$-th) tensor power of itself, given by:
$$
\om_d: \la \longrightarrow {\la}^{\otimes d}
\eqno(5.1)
$$
where $\om_d$ on any fiber of $\la$ is the map $z \mapsto z^{d}$.
Observe that $\om_d$ maps $\la$ minus its zero section to
${\la}^{\otimes d}$ minus its zero section by a map which
is of degree $d$ on the $\CC^{*}$ fibers.
{\it $\om_d$ is a homomorphism of the associated principal $\CC^{*}$
bundles}. When $T$ is a complex manifold, and $\la$ is
a line bundle in that category, then the map $\om_d$ is a holomorphic
morphism between the total spaces of the source and target bundles.

In the situation of Lemma V.1, therefore, we may define a canonical
principal bundle morphism relating the relevant bundles:
$$
\Om(\pi) := \Ga(\pi) \circ \om_{deg(\pi)}:\Dng \rightarrow \Dnt
\eqno(5.2)
$$
where $\Ga(\pi)$ is the canonical (GRR) line bundle morphism provided by
the Lemma.

The canonical and functorial choice of these connecting maps,
$\Om(\pi)$, provides us with a direct system of line/principal
bundles over the direct systems of Teichm\"uller spaces.

Given a direct system $T_{\a}$ of complex manifolds, and line
bundles $\xi_{\a}$ over these, whenever there are connecting
maps as the $\Om(\pi)$ above, we may pass to the direct limit of the
bundles themselves, simultaneously with passing to the limit of the
base spaces. That precipitates our main ``non-perturbative'' result:

\noindent
{\bf Theorem V.2:} {\it Fix integer $n$. Starting from any
``base'' surface $X \in Ob(\A)$, we obtain a direct system of principal
$\CC^*$ bundles $\L_n(\tX) := DET_{n,g(\tX)}$ over the
Teichm\"uller spaces $\T(\tX)$ with connecting holomorphic
homomorphisms $\Om(\pi)$ between their total spaces.

Passing to the direct limit, one therefore obtains over the
universal commensurability Teichm\"uller space,
$\TinX$, a principal $\CQ$ bundle:
$$
\L_{n,\infty}(X) = ind. lim. \L_n(Y)
$$
Since the maps $\Om(\pi)$ preserved the Quillen unit circles, the limit
object also inherits such a Quillen ``hermitian'' structure.
The commensurability modular group action $CM_{\ify}(X)$ on
$\TinX$ has a natural lifting to $\L_{n,\infty}(X)$ -- acting by isometrical
automorphisms.

Finally, the Mumford isomorphisms persist as isometrical $\MCin$
equivariant isomorphisms:
$$
\L_{n,\infty}(X) = (6n^{2} - 6n + 1) \L_{1,\infty}(X)
$$}

\noindent
{\it Remarks:} (i) To be explicit, the Mumford
isomorphism in the above theorem means that
$\L_{n,\infty}$ and $\L_{1,\infty}$ are equivariantly isomorphic
relative to the automorphism of $\CQ$ induced by the homomorphism
of $\C*$ that raises to the power exhibited.

(ii) We could have used the Quillen hermitian structure to reduce the
structure group from $\CC^{*}$ to $U(1)$, and thus obtain direct systems
of $U(1)$ bundles over the Teichm\"uller spaces. Passing to the direct
limit would then produce $U(1) \otimes \QQ := $ ``tiny circle'' bundles
over $\Tin$.

(iii) Another interpretation of this construction over $\Tin$ is to
produce ``rational line bundles'' over it, with relating Mumford
isomorphisms, utilizing Lemma V.1 but not involving the power law
mappings. See [BNS] for details.

(iv) We have also glued together the universal bundles over direct limits
of the {\it moduli} spaces $\Mg$, by considering the tower of {\it
characteristic} coverings over $X$. See [BN].

Above we have succeeded in fitting together the Hodge and higher $DET_n$
bundles over the ind-space $\Tin$, together with the relating Mumford
isomorphisms -- our entire construction being $\MCin$ equivariant.
We thus have a structure on $\Tin$ that suggests a {\it genus-independent
and universal} version of the finite dimensional Polyakov structure that
we delineated in Sections I to III of this paper.


\noindent
{\bf Author's e-mail: nag@imsc.ernet.in}


\baselineskip12pt

\begin{thebibliography}{99}

\bibitem[Alv]{Alv} O. Alvarez, Theory of strings with boundaries:
fluctuations, topology and quantum geometry, {\it Nucl. Phys.
B216}, (1983), 125-184.

\bibitem[AN]{AN} L. Alvarez-Gaume and P. Nelson, Riemann surfaces
and string theories, in {\it Trieste Spring School ``Supersymmetry,
supergravity,superstrings '86''} World Scient., Singapore, (1986),
419-510.

\bibitem[BGS]{BGS} J.Bismut, H.Gillet and C.Soul\'e, Analytic
torsion and holomorphic determinant bundles, I,II,III, {\it
Comm. Math Phys., 115}, 1988, 49-126, 301-351.

\bibitem[BK]{BK} A. Belavin, V. Knizhnik, Complex geometry and
quantum string theory. {\it Phys. Lett., 168 B}, (1986),
201-206.

\bibitem[BN]{BN} I. Biswas and S. Nag, Weil-Petersson geometry
and determinant bundles over inductive limits of moduli spaces,
{\it preprint}.

\bibitem[BNS]{BNS} I. Biswas, S. Nag  and D. Sullivan, Determinant
bundles, Quillen metrics and Mumford isomorphisms over the Universal
Commensurability Teichm\"uller Space, {\it IHES (Paris) preprint,
IHES/M/95/43}, May 1995.

\bibitem[Bos]{Bos} J.B. Bost, Fibres determinants, determinants
regularises et mesures sur les espaces de modules des courbes
complexes, {\it Semin.  Bourbaki, 152-153}, (1987), 113-149.

\bibitem[D]{D} P. Deligne, Le d\'eterminant de la cohomologie,
{\it Contemporary Math., 67}, (1987), 93-177.

\bibitem[DP]{DP} E. D'Hoker and D. Phong, Multiloop amplitudes for the
bosonic Polyakov string, {\it Nucl Phys B, 269}, (1986), 205-234.

\bibitem[IT]{IT} Y. Imayoshi and M. Taniguchi, {\it An introduction to
Teichm\"uller spaces}, Springer-Verlag, Tokyo, (1992).

\bibitem[J]{J} J. Jost, Harmonic maps and curvature computations in
Teichm\"uller theory, {\it Ann. Acad. Scient Fenn., 16}, (1991)), 13-46.

\bibitem[KM]{KM} F. Knudsen and D. Mumford, The projectivity of
the moduli space of stable curves I, preliminaries on ``det''
and ``div'', {\it Math. Scand. 39}, (1976) 19-55.

\bibitem[Ma]{Ma} A. Mann, Problem about automorphisms of infinite groups,
{\it Second Intn'l Conf Group Theory}, Debrecen, (1987).

\bibitem[MT]{MT} F. Menegazzo and M. Tomkinson, Groups with trivial
virtual automorphism group, {\it Israel Jour. Math, 71}, (1990), 297-308.

\bibitem[Mum]{Mum} D. Mumford, Stability of projective varieties,
{\it L'Enseign. Math., 23}, (1977), 39-100.

\bibitem[N1]{N1} S. Nag, {\it The Complex Analytic Theory of Teichm\"uller
Spaces},  Wiley-Interscience, New York, (1988).

\bibitem[N2]{N2} S. Nag, Canonical measures on the moduli spaces of
compact Riemann surfaces, {\it Proc. Indian Acad. Sci., (Math Sci.), 99}
(1989), 103-111.

\bibitem[N3]{N3} S. Nag, Self-dual  Connections,  hyperbolic  metrics
and harmonic Mappings  on   Riemann surfaces, {\it Proc. Indian Acad. Sci.,
(Math Sci.), 101}, (1991), 215-218.

\bibitem[N4]{N4} S. Nag, A period mapping in Universal Teichm\"uller
space, {\it Bull. Amer. Math. Soc., 26, no.2}, (1992), 280-287.

\bibitem[NS]{NS} S. Nag, D. Sullivan, Teichm\"uller theory and
the universal period mapping via quantum calculus and the
$H^{1/2}$ space on the circle, {\it Osaka J. Math., 32, no.1},
(1995), 1-34.  [Max-Planck-Inst. (Bonn) preprint MPI-94-54.]

\bibitem[NV]{NV} S. Nag and A. Verjovsky,
 Diff$(S^{1})$ and the Teichm\"uller Spaces, (Parts I and  II),
{\it Commun. Math. Phys., 130}, (1990), 123-138. Part I with A.V.

\bibitem[Nel]{Nel} P.Nelson, Lectures on strings and moduli
space, {\it Physics Reports, 149, no. 6}, (1987), 337-375.

\bibitem[P]{P} R.Penner, Universal constructions in Teichm\"uller theory,
{\it Adv. Math. 98}, (1993), 143-215.

\bibitem[Pol1]{Pol 1} A.M. Polyakov, Quantum geometry of bosonic
strings, {\it Physics Letters, 103B, no. 3}, (1981), 207-210.

\bibitem[Pol2]{Pol2} A.M. Polyakov, Quantum geometry of fermionic
strings, {\it Physics Letters, 103B, no. 3}, (1981), 211-213.

\bibitem[Q]{Q} D. Quillen, Determinants of Cauchy-Riemann
operators over a Riemann surface. {\it Func. Anal. Appl., 19},
(1985), 31-34.

\bibitem[S]{S} D. Sullivan, Relating the universalities of
Milnor-Thurston, Feigenbaum and Ahlfors-Bers, in Milnor
Festschrift {\it Topological Methods in Modern mathematics} (ed.
L. Goldberg, A. Phillips) Publish or Perish, (1993), 543-563.

\bibitem[Sha]{Sha} I.R. Shafarevich, On some
infinite-dimensional groups II, {\it Math USSR Izvest., 18},
(1982), 185-194.

\bibitem[W]{W} M. Wolf, The Teichm\"uller theory of harmonic maps,
{\it Jour. Diff. Geom., 29}, (1989), 449-479.

\bibitem[Wol]{Wol} S. Wolpert, Chern forms and the Riemann
tensor for the moduli space of curves, {\it Invent. Math., 85},
(1986), 119-145.

\bibitem[ZT]{ZT} P.G. Zograf and L.A. Takhtadzhyan, A local
index theorem for families of {$\bar \partial$}- operators on
Riemann surfaces, {\it Russian Math Surveys, 42}, (1987),
169-190.

\end{thebibliography}
\end{document}
